\documentclass[12pt]{article}
\usepackage{latexsym}
\usepackage{amssymb}
\usepackage{amsthm}
\usepackage{amsmath}
\usepackage{color}
\usepackage{url}
\usepackage{comment}
\usepackage{graphicx}

\includecomment{jp-text}
\excludecomment{en-text}

\newif\ifmemo
\memofalse

\newif\ifrs
\rstrue
\ifrs \usepackage{mathrsfs} \fi

\newif\ifcol
\coltrue
\def\E{{\bf E}}
\def\V{{{\bf Var}}}
\def\de{{\rm d}}
\def\I{{{\bf 1}}}

\newcommand{\colorr}{\color[rgb]{0,0,0}}

\newcommand{\colorb}{\color[rgb]{0,0,0}}
\newcommand{\colord}{\color[rgb]{0,0,0}}
\newcommand{\colorn}{\color[rgb]{1,1,1}}

\begin{document}
%nakamacro.tex(H120522;0730)
%\documentstyle[11pt]{article}
%\setlength{\textwidth}{6.5in}
%\setlength{\oddsidemargin}{0in}
%\setlength{\topmargin}{-0.52in}
%\setlength{\textheight}{9.0in}
%\setlength{\footskip}{0.7in}

%\input nakamacro-s.tex
%%%%%%%%%%%%%%%%%%%%%%%%%%%%%%%%%%%%%%%%%%%%%%%%%%%%%%%%%%%%%%%%%%%
%%%%%%%%%%%%%%%%%%%%%%%%%%%%%%%%%%%%%%%%%%%%%%%%%%%%%%%%%%%%%%%%%%%
%%%%%%%%%%%%%%%%%%%%%%%%%%%%%%%%%%%%%%%%%%%%%%%%%%%%%%%%%%%%%%%%%%%

% == Bold Math ==
%\newcommand{\mb}[1]{\mbox{\boldmath $ #1 $}}
\newcommand{\mb}[1]{{\bf #1 }}

% == Black Board ==
% == Black Board ==
\newcommand{\bbA}{{\mathbb A}}
\newcommand{\bbB}{{\mathbb B}}
\newcommand{\bbC}{{\mathbb C}}
\newcommand{\bbD}{{\mathbb D}}
\newcommand{\bbE}{{\mathbb E}}
\newcommand{\bbF}{{\mathbb F}}
\newcommand{\bbG}{{\mathbb G}}
\newcommand{\bbH}{{\mathbb H}}
\newcommand{\bbI}{{\mathbb I}}
\newcommand{\bbJ}{{\mathbb J}}
\newcommand{\bbK}{{\mathbb K}}
\newcommand{\bbL}{{\mathbb L}}
\newcommand{\bbM}{{\mathbb M}}
\newcommand{\bbN}{{\mathbb N}}
\newcommand{\bbO}{{\mathbb O}}
\newcommand{\bbP}{{\mathbb P}}
\newcommand{\bbQ}{{\mathbb Q}}
\newcommand{\bbR}{{\mathbb R}}
\newcommand{\bbS}{{\mathbb S}}
\newcommand{\bbT}{{\mathbb T}}
\newcommand{\bbU}{{\mathbb U}}
\newcommand{\bbV}{{\mathbb V}}
\newcommand{\bbW}{{\mathbb W}}
\newcommand{\bbX}{{\mathbb X}}
\newcommand{\bbY}{{\mathbb Y}}
\newcommand{\bbZ}{{\mathbb Z}}

%\newcommand{\bbd}{{\mathbb d}}?

%%%%%%%%%%%%%%%%%%%%%%%%%%%%%%%%%%%%%%%%%%%%%%%%%%%%%%%%%%%%%%%%%%%
%%%%%%%%%%%%%%%%%%%%%%%%%%%%%%%%%%%%%%%%%%%%%%%%%%%%%%%%%%%%%%%%%%%
%%%%%%%%%%%%%%%%%%%%%%%%%%%%%%%%%%%%%%%%%%%%%%%%%%%%%%%%%%%%%%%%%%%

\newtheorem{definition}{Definition}
\newtheorem{assumption}{$[$ A}
\newtheorem{condition}{$[$ C}
\newtheorem{lemma}{Lemma}
\newtheorem{proposition}{Proposition}
\newtheorem{theorem}{Theorem}
\newtheorem{remark}{Remark}
\newtheorem{example}{Example}
\newtheorem{cor}{Corollary}
\newtheorem{general}{}
%--------------------------------------------------------------------------
%BOLD FACES
\def\bsx{\bar{x}}
\def\bsd{\bar{d}}
\def\bx{\bar{X}}
\def\ba{\bar{A}}
\def\bb{\bar{B}}
\def\bc{\bar{C}}
\def\bv{\bar{V}}
\def\n{{\bf n}}
\def\A{{\bf A}}
\def\B{{\bf B}}
\def\C{{\bf C}}
\def\D{{\bf D}}
\def\E{{\bf E}}
\def\F{{\bf F}}
\def\G{{\bf G}}
\def\I{{\bf I}}
\def\J{{\bf J}}
\def\K{{\bf K}}
\def\L{{\bf L}}
\def\M{{\bf M}}
\def\N{{\bf N}}
\def\O{{\bf O}}
\def\P{{\bf P}}
\def\Q{{\bf Q}}
\def\R{{\bf R}}
\def\S{{\bf S}}
\def\T{{\bf T}}
\def\U{{\bf U}}
\def\V{{\bf V}}
\def\W{{\bf W}}
\def\X{{\bf X}}
\def\Y{{\bf Y}}
\def\Z{{\bf Z}}
\def\cala{{\cal A}}
\def\calb{{\cal B}}
\def\calc{{\cal C}}
\def\cald{{\cal D}}
\def\cale{{\cal E}}
\def\calf{{\cal F}}
\def\calg{{\cal G}}
\def\calh{{\cal H}}
\def\cali{{\cal I}}
\def\calj{{\cal J}}
\def\calk{{\cal K}}
\def\call{{\cal L}}
\def\calm{{\cal M}}
\def\caln{{\cal N}}
\def\calo{{\cal O}}
\def\calp{{\cal P}}
\def\calq{{\cal Q}}
\def\calr{{\cal R}}
\def\cals{{\cal S}}
\def\calt{{\cal T}}
\def\calu{{\cal U}}
\def\calv{{\cal V}}
\def\calw{{\cal W}}
\def\calx{{\cal X}}
\def\caly{{\cal Y}}
\def\calz{{\cal Z}}
%
%YOKUTUKAUMONO
\def\hd{\hat{\bbD}}
\def\cosech{\mbox{\rm \,cosech}}
\def\sech{\mbox{\rm \,sech}}
\def\sskip{\hspace{0.5cm}}
\def\simleq{ \raisebox{-.7ex}{\em $\stackrel{{\textstyle <}}{\sim}$} }
\def\leqsim{ \raisebox{-.7ex}{\em $\stackrel{{\textstyle <}}{\sim}$} }
\def\ep{\epsilon}
\def\half{\frac{1}{2}}
\def\iku{\rightarrow}
\def\Iku{\Rightarrow}
\def\ikup{\rightarrow^{p}}
\def\inclusion{\hookrightarrow}
\def\cadlag{c\`adl\`ag\ }
\def\up{\uparrow}
\def\down{\downarrow}
\def\doti{\Leftrightarrow}
\def\douti{\Leftrightarrow}
\def\dochi{\Leftrightarrow}
\def\douchi{\Leftrightarrow}%
%KAIGYOU,ARRAY
\def\yy{\\ && \nonumber \\}
\def\y{\vspace*{3mm}\\}
\def\nn{\nonumber}
\def\be{\begin{equation}}
\def\ee{\end{equation}}
\def\bea{\begin{eqnarray}}
\def\eea{\end{eqnarray}}
\def\beas{\begin{eqnarray*}}
\def\eeas{\end{eqnarray*}}
%
%KONO RONBUN DE TUKAU MONO
\def\]]{]]}
\def\[[{[[}
\def\sskip{\hspace{5mm}}
\def\proof{{\it Proof. }}
\def\balpha{\bar{\alpha}}
\def\bbalpha{\bar{\bar{\alpha}}}
\def\combi{\l(\begin{array}{c}\alpha\\ \beta \end{array}\r)}
\def\f{^{(1)}}
\def\s{^{(2)}}
\def\ss{^{(2)*}}
\def\l{\left}
\def\r{\right}
\def\a{\alpha}
\def\b{\beta}
\def\limsup{\overline{\lim}}
\def\liminf{\underline{\lim}}
%上に定義されたコマンドは数式モ−ドで用いる。
%--------------------------------------------------

%\setcounter{page}{0}

%\begin{itemize}
%\item {\colorr Correction/Question}
%\item {\colorb Proposal}
%\item {\colorg  Believed without checking. Please check once again.} 
%\item We should check the validity of  the original paper we referred on the uniform HR theroem. 
%\end{itemize}

%\newpage 

\title{Estimation for the change point of the volatility 
in a stochastic differential equation
\footnote{
This work was 
in part supported by 
JST Basic Research Programs PRESTO, 
Grants-in-Aid for Scientific Research 
No. 19340021, 
the Global COE program 
``The research and training center for new development in mathematics'' 
of 
Graduate School of Mathematical Sciences, University of Tokyo, 
and by Cooperative Research Program 
of the Institute of Statistical Mathematics.} 
}

\author{Stefano M. Iacus\footnote{Corresponding author.}\\
{\small Department of Economics, Business and Statistics, University of Milan}\\ 
{\small Via Conservatorio 7, 20122 Milan, Italy; 
  \url{stefano.iacus@unimi.it}}
\and
Nakahiro Yoshida\\
{\small University of Tokyo, 
and 
Japan Science and Technology Agency}\\
{\small Graduate School of Mathematical Sciences, University of Tokyo}\\
{\small 3-8-1 Komaba, Meguro-ku, Tokyo 153-8914 Japan; \url{nakahiro@ms.u-tokyo.ac.jp}}
}


\maketitle

\vspace{-1cm}
\begin{abstract}
We consider a multidimensional It\^o process $Y=(Y_t)_{t\in[0,T]}$ 
with some unknown drift coefficient process $b_t$ and 
volatility coefficient $\sigma(X_t,\theta)$ 
with covariate process $X=(X_t)_{t\in[0,T]}$, 
the function $\sigma(x,\theta)$ being 
known up to $\theta\in\Theta$. 
For this model we consider a change point problem for the parameter $\theta$ 
in the volatility component. 
The change is supposed to occur at some point $t^*\in (0,T)$. 
Given discrete time observations from the process $(X,Y)$, 
we propose quasi-maximum likelihood estimation of the change point. 
We present the rate of convergence of the change point estimator 
and the limit thereoms of aymptotically mixed %normal 
type. 
\end{abstract}
{\bf keywords} It\^o processes, discrete time observations, change point estimation, volatility

\section{Introduction}
The problem of change point has been considered initially in the framework of independent and identically distributed data  by many authors, see e.g. Hinkley (1971), Cs\"{o}rg\H{o} and Horv\'{a}th (1997), Inclan and Tiao (1994).
Recently, it naturally moved to context of time series analysis, see for example, Kim {\it et al.} (2000), Lee {\it et al.} (2003), Chen {\it et al.} (2005) and the papers cited therein.

In fact, change point problems have originally arisen in the context of quality control, but the problem of abrupt changes in general arises in many contexts like epidemiology, rhythm analysis in electrocardiograms, seismic signal processing, study of archeological sites and financial markets. In particular, in the analysis of financial time series, the knowledge of the change in the volatility structure of the process under consideration is of a certain interest.

In this paper we deal with a change-point problem for the volatility 
of a process solution to a stochastic differential equation, 
when observations are collected at discrete times. 
The instant of the change in volatility regime is identified retrospectively 
by maximum likelihood method on the approximated likelihood. 
For continuous time observations of 
diffusion processes Lee {\it et al.} (2006) 
considered the change point estimation problem for the drift.  
In the present work we only assume regularity conditions on the drift process. 
De Gregorio and Iacus (2008) considered a least squares approach following 
the lines of Bai (1994, 1997) of a simplified model also under discrete sampling {\colorb while Song and Lee (2009) considered a CUSUM approach. }
Finally it should be noted that the problems of the change-point of 
drift for ergodic diffusion processes have been treated by Kutoyants (1994, 2004), 
but the asymptotics and 
the sampling schemes are different from this paper. 

The paper is organized as follows. Section \ref{sec:model} introduces the model of observation, the regularity conditions and some notation.
Section \ref{sec:consist} studies consistency and the rate of convergence of estimator of the change while asymptotic distributions are considered 
in Section \ref{sec:asym}. 
{\colord{
A mixture of certain Wiener functionals appears as the limit of 
the likelihood ratio random field, and it characterizes 
the limit distribution of the change-point estimator. }}
Those sections assume that consistent estimators of the volatility parameters  are available. Section \ref{081227-1} presents 
some practical considerations and a proposal to obtain first stage estimators of the volatility parameters which allow to obtain all asymptotic properties stated in the previous sections.
Finally, Section \ref{201229-1} presents some 
numerical analysis to asses the performance of the estimators. Tables are collected at the end of the paper.

\section{Estimator for the change-point of the volatility}\label{sec:model}
%\section{Model and assumptions}\label{sec:model}
%Given a stochastic basis $\calb=(\Omega,\calf,\F,P)$ 
%with filtration $\F=(\calf_t)_{t\in[0,T]} 
Consider a {\colord{$d$-}}dimensional It\^o process described by the 
stochastic differential equation 
%It\^o process $\{X_t, 0\leq t \leq T\}$ $X=(X_t)_{t\in[0,T]}$ given by 
\begin{equation}
d Y_t = b_t d t + \sigma(X_t,\theta) d W_t,\ \ t\in[0,T], 
\label{eq1}
\end{equation}
where $W_t$ is {\colord{an $r$-dimensional}} standard Wiener process, 
on a stochastic basis, 
$b_t$ and $X_t$ are {\colord{vector valued}} 
progressively measurable processes, 
and $\sigma(x,\theta)$ is a {\colord{matrix valued}} 
function. %given $\theta$ . 

We assume that 
there is the time $t^*$ across which the diffusion coefficient 
changes from $\sigma(x,\theta_0)$ to $\sigma(x,\theta_1)$. 
The change point $t^*\in(0,T)$ is unknown and 
we want to estimate $t^*$ based on the observations sampled 
from the path of $(X,Y)$. 
The coefficient $\sigma(x,\theta)$ is assumed to be known 
up to the parameter $\theta$, while 
$b_t$ is completely unknown and unobservable, therefore possibly 
depending on $\theta$ and $t^*$.

The sample consists of $(X_{t_i},Y_{t_i})$, $i=0,1,...,n$, where 
$t_i=ih$ for $h=h_n=T/n$. 
The parameter space $\Theta$ of $\theta$ 
is a bounded domain in $\bbR^{d_0}$, $d_0\geq 1$, 
{\colord{and}} the parameter 
$\theta$ is a nuisance in estimation of $t^*$. 
Denote by $\theta_i^*$ the true value of $\theta_i$ 
for $i=0,1$. 

Let $\vartheta_n=|\theta_1^*-\theta_0^*|$. 
We will consider the following two different situations. 
\begin{description}
\item[(A)] 
$\theta_0^*$ and $\theta_1^*$ 
are fixed and do not depend on $n$. 

\item[(B)] 
$\theta_0^*$ and $\theta_1^*$ depend on $n$, and as $n\to\infty$, 
$\theta_0^*\to\theta^*\in\Theta$, 
$\vartheta_n\to0$ and 
$n\vartheta_n^2\to\infty$.% for some constant $\gamma>1$. 
\end{description}
In Case {\bf (A)}, $\vartheta_n$ is a constant $\vartheta_0$ 
independent of $n$.

{\colord{We shall formulate the problem more precisely.}} 
It will be assumed that 
the process $Y$ generating the data 
is an It\^o process realized on a stochastic basis 
$\calb=(\Omega,\calf,\F,P)$ with filtration $\F=(\calf_t)_{t\in[0,T]}$, 
and satisfies the stochastic integral equation 
\beas 
Y_t&=&
\Bigg\{
\begin{array}{ll}
Y_0+\int_0^t b_s d s+\int_0^t \sigma(X_s,\theta_0^*) d W_s
& \mbox{ for } t\in[0,t^*)
\y
Y_{t^*}+\int_{t^*}^t b_s d s+\int_{t^*}^t \sigma(X_s,\theta_1^*) d W_s
& \mbox{ for } t\in[{t^*},T]. 
\end{array}
\eeas
Here $W_t$ is an $r$-dimensional $\F$-Wiener process on $\calb$, 
and $b_t$, $X_t$ and $\sigma(x,\theta)$ satisfy the conditions below. 
Let $\calx$ be a closed set in $\bbR^{d_1}$ 
{\colord{(possibly $\calx=\bbR^{d_1}$)}} 
and denote the modulus of continuity of a function 
$f:I\to\bbR^{d_1}$ by 
\beas 
w_I(\delta,f)=\sup_{s,t\in I,\>|s-t|\leq\delta}
|f(s)-f(t)|.
\eeas 
For matrices $A=(a_{ij})$ and $B=(b_{ij})$ of the same size, 
we write $A^{\otimes2}=A\>{^{\sf t}A}$, 
$A[B]=\sum_{ij}a_{ij}b_{ij}=\mbox{Tr}(A\>{^{\sf t}B})$, 
and the Euclidean norm of $A$ by $|A|=(A[A])^{1/2}$. 
Set $S(x,\theta)=\sigma(x,\theta)^{\otimes2}$. 
$\lambda_1(A)$ denotes the minimum eigenvalue of 
a symmetric matrix $A$.

\begin{description}
\item[[H\!\!]]$_j$ 
(i) $\sigma(x,t)$ is a 
measurable function defined on $\calx\times[0,T]$ satisfying 
\begin{itemize}
\item[(a)] 
$\inf_{(x,\theta)\in\calx\times\Theta}\lambda_1(S(x,\theta))>0$, 
\item[(b)] derivatives 
$\partial_\theta^\ell\sigma$ ($0\leq\ell\leq j {\colorr +[d_0/2]}$) 
exist and those functions are continuous 
on $\calx\times\Theta$, 
\item[(c)] there exists a locally bounded function 
$L:\calx\times \calx\times \Theta\to\bbR_+$ such that 
\beas 
|\sigma(x,\theta)-\sigma(x',\theta)|
\leq L(x,x',\theta)|x-x'|^\alpha\sskip(x,x'\in \calx,\>\theta\in\Theta) 
\eeas
for some constant $\alpha>0$. %\geq(\gamma-1)/\beta$. 
\end{itemize}
%
%\vspace{1mm}\\ 
\noindent 
(ii) $(X_t)_{t\in[0,T]}$ is a %$d$-dimensional 
progressively measurable process 
taking values in $\calx$ such that 
%\beas 
%\sup_{s,t\in[0,T]:|s-t|\leq n^{-1}}|X_s-X_t|=o_p(\vartheta_n^\beta) 
%\eeas
\beas 
w_{[0,T]}\bigg(\frac{1}{n},X\bigg)=o_p(\vartheta_n^{1/\alpha})
\eeas
as $n\to\infty$. 
\vspace{1mm}\\ \noindent
(iii) $(b_t)_{t\in[0,T]}$ is a progressively measurable process 
taking values in $\bbR^d$ 
such that 
$(b_t-b_0)_{t\in[0,T]}$ is locally bounded. \\

\begin{en-text}
$[$ {\bf We will add condition later, but the localization technique 
will be applied. }$]$ 
\end{en-text}
\end{description}

\begin{en-text}
The symbol $d_s(\calg)$ stands for the $\calg$-stable convergence 
in distribution (see e.g. Jacod, 1997). 
\end{en-text}

\begin{remark}\rm 
The term ``locally bounded'' in [H]$_j$ (i) (c) means, as usual, 
being bounded on every compact set. 
The case where the drift $b_t$ changes its structure 
at time $t^*$, or any time in force, 
is included in our context because 
$b_t$ admits jumps. 
The case of time dependent $\sigma$ is included by making 
$X_t$ have argument $t$. 
Needless to say, if we set $X$ or a part of $X$ as $Y$, then 
our model can express a system with feedback, 
{\colord{in particular, a diffusion process}}. 
By [H]$_j$ (ii), $t\mapsto X_t$ is continuous a.s. 
Also, [H]$_j$ (ii) imposes a restriction on the rate $\vartheta_n$. 
For example, when $\alpha=1$, for a Brownian motion $X$, 
it suffices 
that $n\vartheta_n^2/\log n\to\infty$, due to L\'evy property. 
The additional ${\colorr [d_0/2]}$ time differentiability 
to $j$ is used only in Step (iii) of the proof of Theorem \ref{081219-1}. 
Therefore, it is possible to replace the range of $\ell$ 
to ``$0\leq\ell\leq j$'' under a condition that ensures 
the the H\'ajek-Renyi type estimate 
just before going to 
Inequality (\ref{201221-11}) below.
\end{remark}

\begin{en-text}
We assume that 

when $\theta_0^*$ and $\theta_1^*$ are fixed, that is, independent of $n$. 
Those parameters are possibly unknown. 
It should be noted that consistent estimation is impossible 
without a condition like the last one.

In case $\vartheta_n\iku0$ and the parameters are known, 
$h(s)$ is assumed to be 
a continuous process satisfying 
$h(t^*)>0$ a.s. and 
\bea\label{200124-1} 
0\leq h(s) &\leq&
\liminf_{n\to\infty} 
\vartheta_n^{-2}\left[ 
\frac{\sigma(s;\theta_1^*)^2}{\sigma(s;\theta_0^*)^2} - 1 
-\log
\frac{\sigma(s;\theta_1^*)^2}{\sigma(s;\theta_0^*)^2}
\right]\sskip a.s. 
\eea

In case $\vartheta_n\iku0$ and the parameters are unknown, 
$h(s)$ is assumed to be 
a continuous process satisfying 
$h(t^*)>0$ a.s. and 
\bea\label{200130-2} 
0\leq h(s) &\leq&
2\inf_{\theta}
\left( 
\frac{\partial_\theta \sigma(s;\theta)}{\sigma(s;\theta)}
\r)^2 \sskip a.s. 
\eea
If $\theta_0$ is bounded, it is possible to restrict 
the set where infimum is taken to 
a neibourghood of the set where $\theta_0$ runs. 
\end{en-text}

\begin{en-text}
We assume that
\beas 
P\l[\sigma(t^*;\theta_0^*)\not=\sigma(t^*;\theta_1^*)\r]=1
\eeas

We assume that the existence of the continuous process
$h(s)$  satisfying 
$h(t^*)>0$ a.s. and 
\bea\label{200124-1} 
0\leq h(s) &\leq&
\liminf_{n\to\infty} 
\vartheta_n^{-2}\left[ 
\frac{\sigma(s;\theta_1^*)^2}{\sigma(s;\theta_0^*)^2} - 1 
-\log
\frac{\sigma(s;\theta_1^*)^2}{\sigma(s;\theta_0^*)^2}
\right]\sskip a.s. 
\eea
\end{en-text}

%\section{Estimation}\label{sec:est}

%{\bf koko} 

Write $\Delta_iY=Y_{t_i}-Y_{t_{i-1}}$ and 
let 
\beas 
\Phi_n(t;\theta_0,\theta_1)
&=&
\sum_{i=1}^{[nt{\colorr /T}]}G_i(\theta_0)+\sum_{i={[nt{\colorr /T}]+1}}^nG_i(\theta_1),
\eeas
where 
\beas 
G_i(\theta) &=&
\log\det S(X_{t_{i-1}},\theta)
+ h^{-1}S(X_{t_{i-1}},\theta)^{-1}[(\Delta_iY)^{\otimes2}]. 
\eeas

Suppose that there exists an estimator 
$\hat{\theta}_k$ for {\colord{each}} $\theta_k$, $k=0,1$. 
Each estimator is based on $(X_{t_i},Y_{t_i})_{i=0,1,...,n}$ and so 
depends on $n$. 
To make our discussion complete, 
in case $\theta_k^*$ are known, we define $\hat{\theta}_k$ 
just as $\hat{\theta}_k=\theta_k^*$. 
This article proposes 
\beas 
\hat{t}_n&=&\mbox{argmin}_{t\in[0,T]} 
\Phi_n(t;\hat{\theta}_0,\hat{\theta}_1)
\eeas
for the estimation of $t^*$. 
More precisely, $\hat{t}_n$ is 
any measurable function of $(X_{t_i})_{i=0,1,...,n}$ 
satisfying 
\beas 
\Phi_n(\hat{t}_n;\hat{\theta}_0,\hat{\theta}_1)
&=&
\min_{t\in[0,T]} 
\Phi_n(t;\hat{\theta}_0,\hat{\theta}_1). 
\eeas

%%%%%%%%%%%%%%%%%%%%%%%%%%%%%%%%%%%%%%%%%%%%%%%%%%%%%%%%%%%%
%%%%%%%%%%%%%%%%%%%%%%%%%%%%%%%%%%%%%%%%%%%%%%%%%%%%%%%%%%%%
%%%%%%%%%%%%%%%%%%%%%%%%%%%%%%%%%%%%%%%%%%%%%%%%%%%%%%%%%%%%
\section{Rate of convergence}\label{sec:consist}
We introduce identifiability conditions in order to ensure  consistent 
estimation. In Case {\bf (A)} we assume

\begin{description}
\item[[A\!\!]] %$_{\bf I}$ 
%$\sigma(t;\theta)$ are positive and continuous, and 
$
P\l[S(X_{t^*};\theta_0^*)\not=S(X_{t^*};\theta_1^*)\r]=1
$; 
\end{description}
In Case {\bf (B)} we assume 
%{\bf 条件検討！}
%
\begin{description}
\item[[B\!\!]] %$_{\bf II}$
$\Xi(X_{t^*},\theta^*)$ is positive-definite a.s., 
where 
\beas 
\Xi(x,\theta)=
\bigg(\mbox{Tr}((\partial_{\theta^{(i_1)}} S)S^{-1}
(\partial_{\theta^{(i_2)}} S)S^{-1})(x,\theta)
\bigg)_{i_1,i_2=1}^{d_0}, \sskip
\theta=(\theta^{(i)}). 
\eeas 
\end{description}

\begin{remark}\rm 
Since $\Xi(x,\theta^*)$ is the Hessian matrix of 
the nonnegative function 
\beas 
Q(x,\theta^*,\theta):=
\mbox{Tr}\bigg(S(x,\theta^*)^{-1}S(x,\theta)-I_d\bigg)
-\log\det\bigg(S(x,\theta^*)^{-1}S(x,\theta)\bigg)
\eeas
of $\theta$ at $\theta^*$, 
$\Xi(x,\theta^*)$ is nonnegative-definite. 
\end{remark}

\begin{en-text}
\begin{remark}\rm 
Obviously it is sufficient for (\ref{200124-1}) that 
\beas 
0\leq H(x) &\leq&
%\liminf_{n\to\infty} 
|\hat{\theta}_1-\hat{\theta}_0|^{-2}\left[ 
\frac{\sigma(x,\hat{\theta}_1)^2}{\sigma(x,\hat{\theta}_0)^2} - 1 
-\log
\frac{\sigma(x,\hat{\theta}_1)^2}{\sigma(x,\hat{\theta}_0)^2}
\right]\sskip 
\eeas
for all $x\in\bbR$, almost surely, for sufficiently large $n\in\bbN$. 
It is also sufficient for (\ref{200124-1}) that 
\beas
0\leq H(x) &\leq&
|\theta_1-\theta_0|^{-2}\left[ 
\frac{\sigma(x,\theta_1)^2}{\sigma(x,\theta_0)^2} - 1 
-\log
\frac{\sigma(x,\theta_1)^2}{\sigma(x,\theta_0)^2}
\right]\sskip 
\eeas
for all $x\in\bbR$ and $\theta_0,\theta_1\in\Theta$. 
Though these condition have easier expression, 
they are much stronger than (\ref{200124-1}). 

\end{remark}
\end{en-text}

The following property will be necessary 
to validate our estimating procedure. 

\begin{description}
\item[[C\!\!]] 
$
|\hat{\theta}_k-\theta_k^*|=o_p(\vartheta_n)
$ 
as $n\to\infty$ for $k=0,1$. 
\end{description}
In case the parameters are known, 
$\hat{\theta}_k$ should read $\theta_k^*$, and then Condition 
[C] requires nothing. 
Section \ref{081227-1} presents an example of estimator for $\theta_k$ 
which satisfies Condition [C]. 

Here we state the result on the rate of convergence 
of our change-point estimator.

\begin{theorem}\label{081219-1} 
The family 
$\{n\vartheta_n^2(\hat{t}_n-t^*)\}_{n\in\bbN}$ is tight 
under any one of the following conditions. 
\begin{description}
\item[(a)] 
$[H]_1$, $[A]$ and $[C]$ hold in Case {\bf (A)}. 
\item[(b)] 
$[H]_2$, $[B]$ and $[C]$ hold in Case {\bf (B)}. 
\end{description}
\end{theorem}

In Case {\bf (B)}, this result gives consistency of $\hat{t}_n$ 
since $n\vartheta_n^2\to\infty$ by assumption.

The rest of this section will be devoted to the proof of 
Theorem \ref{081219-1}. 
%For this, we may set $T=1$ 
%for notational simplicity without loss of generality. 
Define a stopping time $\tau=\tau(K)$ by 
\beas 
\tau(K) 
&=&
\inf\bigg\{t;\> 
|X_t|+|b_t|%+|\sigma(X_t,\theta_0^*)|+|\sigma(X_t,\theta_1^*)|
>K
\bigg\}\wedge {\colord{T}}
\eeas
for $K>0$. 
$X^\tau$ denotes the process $X$ stopped at $\tau$. 
Write 
$S_i(\theta) = S(X^\tau_{t_i},\theta)$, and 
$\Delta_i Y^\tau = Y^\tau_{t_i} - Y^\tau_{t_{i-1}}$. 
Let 
\beas 
\Psi_n(t;\theta_0,\theta_1)
&=&
\sum_{i=1}^{[nt{\colorr{/T}}]}g_i(\theta_0)+\sum_{i={[nt{\colorr /T}]+1}}^ng_i(\theta_1),
\eeas
where 
\beas 
g_i(\theta) &=&
1_{\{\tau>0\}}
\bigg\{
\log\det S_{i-1}(\theta)
+ h^{-1}S_{i-1}(\theta)^{-1}[(\Delta_iY^\tau)^{\otimes2}]
\bigg\}
\\&=&
1_{\{\tau>0\}}
\log\det S_{i-1}(\theta)
+ h^{-1}S_{i-1}(\theta)^{-1}[(\Delta_iY^\tau)^{\otimes2}]. 
\eeas
Then $\sup_{\theta\in \calk}|g_i(\theta)|\in L^\infty$ 
for any compact set $\calk$ in $\Theta$ under [H]$_1$. 
Denote by $E^{\theta_1^*}_{i-1}$ 
the conditional expectation with respect to 
$\calf_{t_{i-1}}$ 
{\colord{under the true distribution}} 
for $t_{i-1}\geq t^*$.

\begin{lemma}\label{decomposition}
For $t>t^*$, 
\beas 
\Psi_n(t;\theta_0,\theta_1)-\Psi_n(t^*;\theta_0,\theta_1)
&=&
M_n(t;\theta_0,\theta_1)
+A_n(t;\theta_0,\theta_1)
+\rho_n(t;\theta_0,\theta_1),
\eeas
where 
\beas 
M_n(t;\theta_0,\theta_1)
&=&
\sum_{i=[nt^*{\colorr /T}]+1}^{[nt{\colorr /T}]} 
\l\{[g_i(\theta_0)-g_i(\theta_1)]
-E^{\theta_1^*}_{i-1}[g_i(\theta_0)-g_i(\theta_1)]\r\},
\eeas
\beas 
A_n(t;\theta_0,\theta_1)
&=&
1_{\{\tau>0\}}
\sum_{i=[nt^*{\colorr /T}]+1}^{[nt{\colorr /T}]} 
\bigg\{
{\rm Tr}\bigg(S_{i-1}(\theta_0)^{-1}S_{i-1}(\theta_1)-I_d\bigg)
\\&& 
-\log\det\bigg(S_{i-1}(\theta_0)^{-1}S_{i-1}(\theta_1)\bigg) 
\bigg\},
\eeas
\beas
\rho_n(t;\theta_0,\theta_1)
&=&
1_{\{\tau>0\}}
\sum_{i=[nt^*{\colorr /T}]+1}^{[nt{\colorr /T}]}
{\rm Tr} \bigg\{
\l(
S_{i-1}(\theta_1)^{-1}- S_{i-1}(\theta_0)^{-1}\r)\>
\\&&\sskip\sskip
\cdot\l( S_{i-1}(\theta_1)
-h^{-1}E^{\theta_1^*}_{i-1}[(\Delta_iY^\tau)^{\otimes2}]\r)
\bigg\}. 
\eeas
\end{lemma}

The proof of Lemma \ref{decomposition} is straight forward 
and omitted. 

\begin{remark}\rm 
Later we will consider substitution of estimators 
$\hat{\theta}_k$ to $\theta_k$, $k=0,1$. 
Then the expectation 
$E^{\theta_1^*}_{i-1}[g_i(\theta_0)-g_i(\theta_1)]$ is 
taken before the substitution, and so 
\beas 
M_n(t;\hat{\theta}_0,\hat{\theta}_1)
&=&
\sum_{i=[nt^*{\colorr /T}]+1}^{[nt{\colorr /T}]} 
\l\{[g_i(\hat{\theta}_0)-g_i(\hat{\theta}_1)]
-E^{\theta_1^*}_{i-1}[g_i(\theta_0)-g_i(\theta_1)]
\Big|_{\theta_0=\hat{\theta}_0,\theta_1=\hat{\theta}_1}
\r\}. 
\eeas
In particular, the second term in the braces is not necessarily
$\calf_{t_{i-1}}$-measurable. 
\end{remark}

We will need a uniform H\'ajek-R\'enyi inequality. 
Let $D$ be a bounded open set in $\bbR^d$. 
The Sobolev norm is denoted by 
\beas 
\|f\|_{s,p}
&=&
\bigg\{\sum_{i=0}^{s} 
\|\partial_\theta^i f\|_{L^p(D)}^p\bigg\}^{1/p}
\eeas
for $f\in W^{s,p}(D)$, the Sobolev space 
with indices $(s,p)$. 
Suppose that $p>1$ and $s>d/p$. 
The embedding inequality is the following 
\bea\label{201215-1} 
\sup_{\theta\in D} |f(\theta)|
\leq C \|f\|_{s,p}
\sskip(f\in W^{s,p}(D))
\eea
where $C$ is a constant depending only on $s,p$ and $D$. 
We will apply this inequality for $f\in C^s(D)$, 
and the validity of such an inequality 
depends on the regularity of the boundary of $D$; 
see e.g. Yoshida (2005) for the relation to the GRR inequality. 

%uniform HR
\begin{lemma}\label{unif.HR}
Let $(\Omega,\calf,\F=(\calf_j)_{j\in\bbZ_+},P)$ 
be a stochastic basis. 
Let $D$ be a bounded domain in $\bbR^d$ 
admitting Sobolev's inequality {\rm (\ref{201215-1})} 
for some $p\in(1,2]$ and $s\in\bbN$ such that $s>d/p$. 
Let $(c_j)_{j\in\bbZ_+}$ be a nondecreasing sequence of 
positive numbers. 
Let $X=(X_j)_{j\in\bbZ_+}$ 
be a sequence of random fields on $D$ for $j\in\bbZ_+$ 
satisfying the following conditions: 
\begin{description}
\item[(i)] 
For each $(w,j)\in\Omega\times\bbZ_+$, 
$X_j\in C^s(D)$; 
\item[(ii)] 
For each $(\theta,i)\in D\times\{0,1,...,s\}$, 
$(\partial_\theta^i X_j(\theta))_{j\in\bbZ_+}$ 
is a zero-mean $L^p$-martingale with respect to $\F$. 
\end{description}
Then there exists a constant $C'$ depending only on 
$s,p$ and $D$, not depending on $X$, such that 
\beas 
P\bigg[ 
\max_{j\leq n} 
\frac{1}{c_j}\sup_{D}|X_j(\theta)|\geq a \bigg] 
&\leq& 
\frac{C'}{a^p}
\sum_{j=0}^n \frac{1}{c_j^p}
E\bigg[ \|X_j-X_{j-1}\|_{s,p}^p \bigg] 
\eeas
for all $a>0$ and $n\in\bbZ_+$. 
\end{lemma}
\proof 
Let $B=L^p(D)$, then $B$ is $p$-uniformly smooth; 
see Example 2.2 of Woyczy\'nski (1975), p. 247. 
We apply Theorem in Shixin (1997) to conclude 
\beas 
P\bigg[ \max_{j\leq n}\frac{1}{c_j}
\|\partial_\theta^i X_j\|_{{\colord{B}}} \geq a \bigg] 
&\leq&
\frac{C_1}{a^p}
\sum_{j=0}^n \frac{1}{c_j^p} E\bigg[
\|\partial_\theta^i X_j-\partial_\theta^i X_{j-1}\|_{{\colord{B}}}^p \bigg] 
\eeas
for $i\in\{0,1,...,s\}$ for some constant $C_1$. 
Therefore (\ref{201215-1}) yeilds the result. \qed\y

{\it Proof of Theorem \ref{081219-1}.} 
For the proof, we may assume $T=1$ 
for notational simplicity without loss of generality. 
\\\noindent 
(i) Let $\ep$ be an arbitrary positive number. 
Set 
%$H(x)=4c_1\{\sigma(x,\theta_1^*)-\sigma(x,\theta_0^*)\}^2$ 
\beas 
H(x)&=&4Q(x,\theta_0^*,\theta_1^*) 
{\colord{\vartheta_0^{-2}}}
\eeas
in Case {\bf (A)}, 
and set 
$H(x)=\lambda_1(\Xi(x,\theta^*))
%c_1\{\partial_\theta\sigma(x,\theta^*)\}^2
$ in Case {\bf (B)}. 
%The positive constants will be specified later. 
%
We denote 
%$\Psi_n^*(t)=\Psi_n(t;\theta_0^*,\theta_1^*)$, 
$\sigma(t;\theta)=\sigma(X^\tau_t,\theta)$ 
and $h(t)=H(X^\tau_t)$ in what follows. 
Those processes depend on $K$ by definition 
while it is suppressed from the symbols. 
Set $B_K=\{\tau=1\}$ and 
fix a sufficiently large $K$ so that 
$P[B_K^c]<\ep/4$. 

We notice that $h(s)\geq0$ and that 
$h(t^*)>0$ a.s. on $B_K$ 
from the identifiability condition [A]/[B] 
since $X^\tau_{t^*}=X_{t^*}$ on $B_K$. 
We will show that there exists a positive constant $c_\ep$ such that 
\begin{en-text}
\beas
P\l[
\inf_{t\in[t^*,1]}\frac{1}{t-t^*}\int_{t^*}^t 
\l\{
\frac{\sigma(s;\theta_1^*)^2}{\sigma(s;\theta_0^*)^2} - 1 
-\log
\frac{\sigma(s;\theta_1^*)^2}{\sigma(s;\theta_0^*)^2} \r\}\>ds \leq c_\ep
\r]
&<&\ep.
\eeas
\end{en-text}
\beas
P\l[
\inf_{t\in[t^*,1]}\frac{1}{t-t^*}\int_{t^*}^t 
h(s)
\>ds \leq 5c_\ep
\r]
&<&\ep.
\eeas
Define the event $\cala_\delta$ by  
\beas 
\cala_\delta = \left\{ \inf_{t\in[t^*,t^*+\delta]} h(s) 
\geq \half h(t^*) \right\}
\eeas
for $\delta\in(0,1-t^*)$. 
On $\cala_\delta$, it holds that 
\beas 
\inf_{t\in[t^*,t^*+\delta]} \frac{1}{t-t^*}
\int_{t^*}^t h(s) \> ds \geq \half h(t^*)\geq \frac{\delta}{2(1-t^*)}h(t^*)
\eeas
and also that, for $t\in[t^*+\delta,1]$, 
\beas 
\frac{1}{t-t^*}\int_{t^*}^t h(s)\> ds 
&\geq& 
\frac{1}{1-t^*}\int_{t^*}^t h(s)\> ds 
\\&\geq&
\frac{1}{1-t^*}\int_{t^*}^{t^*+\delta} h(s)\> ds 
\\&\geq&
\frac{\delta}{2(1-t^*)}h(t^*).
\eeas

Choose a $\delta$ so that $P[\cala_\delta]>1-\ep/2$ 
by the continuity of $h$, and next 
choose a positive number $c_\ep=c(\ep,\delta)$ such that 
\beas 
P\left[
\frac{\delta}{2(1-t^*)}h(t^*) 
>5c_\ep\right] 
&\geq&
P\left[\bigg\{
\frac{\delta}{2(1-t^*)}h(t^*)>5c_\ep \bigg\}
\bigcap B_K
\right] 
\\&>& 
1-\frac{\ep}{2}. 
\eeas
Then 
\beas 
P\left[\inf_{t\in[t^*,1]}\frac{1}{t-t^*}
\int_{t^*}^t h(s) \>ds\leq 5c_\ep\right]
&\leq&
P[\cala_\delta^c]
+P\left[\cala_\delta,\>
\inf_{t\in[t^*,1]}\frac{1}{t-t^*}
\int_{t^*}^t h(s)\>ds\leq 5c_\ep\right]
\\&<&
\ep.
\eeas

\bigskip \noindent
(ii) With Lemma \ref{decomposition}, 
we decompose 
$\Psi_n(t;\hat{\theta}_0,\hat{\theta}_1)
-\Psi_n(t^*;\hat{\theta}_0,\hat{\theta}_1)$ as follows: 
\beas 
\Psi_n(t;\hat{\theta}_0,\hat{\theta}_1)
-\Psi_n(t^*;\hat{\theta}_0,\hat{\theta}_1)
&=&
M_n(t;\hat{\theta}_0,\hat{\theta}_1)
+A_n(t;\hat{\theta}_0,\hat{\theta}_1)
+\rho_n(t;\hat{\theta}_0,\hat{\theta}_1).
\eeas
Let {\colord{$M\geq1$}}. %$M\in[1,n)$. 
We have 
\bea\label{20080130-1} 
P[n\vartheta_n^2(\hat{t}_n-t^*)>M]
\nn&\leq&
P\l[\inf_{t:n\vartheta_n^2(t-t^*)>M}
\Phi_n(t;\hat{\theta}_0,\hat{\theta}_1)
\leq\Phi_n(t^*;\hat{\theta}_0,\hat{\theta}_1)\r]
\\\nn&\leq&
P\l[\inf_{t:n\vartheta_n^2(t-t^*)>M}
\Psi_n(t;\hat{\theta}_0,\hat{\theta}_1)
\leq\Psi_n(t^*;\hat{\theta}_0,\hat{\theta}_1)\r]
+P[B_K^c]
\\&<&
P_{1,n}+P_{2,n}+P_{3,n}+\ep,
\eea
where
\beas
P_{1,n}&=&
P\l[\sup_{t:n\vartheta_n^2(t-t^*)>M}\frac{1}{[nt]-[nt^*]}
\l|M_n(t;\hat{\theta}_0,\hat{\theta}_1)\r|
\geq\frac{c_\ep\vartheta_n^2}{3}\r]
\\
P_{2,n}&=&
P\l[\inf_{t:n\vartheta_n^2(t-t^*)>M}\frac{1}{[nt]-[nt^*]}
A_n(t;\hat{\theta}_0,\hat{\theta}_1)\leq c_\ep\vartheta_n^2\r]
\\
P_{3,n}&=&
P\l[\sup_{t:n\vartheta_n^2(t-t^*)>M}\frac{1}{[nt]-[nt^*]}
\l|\rho_n(t;\hat{\theta}_0,\hat{\theta}_1)\r|
\geq\frac{c_\ep\vartheta_n^2}{3}\r].
\eeas
Here we read $\inf\emptyset=\infty$ and $\sup\emptyset=-\infty$. 
We will estimate these terms. 

\bigskip\noindent
(iii) Estimate of $P_{1,n}$. 
In Case {\bf (B)}, 
let 
\beas 
\calm_n(t;\theta)
&=&
\sum_{i=[nt^*]+1}^{[nt]} \l\{\partial_\theta g_i(\theta)
-E^{\theta_1^*}_{i-1}[\partial_\theta g_i(\theta)]\r\}.
\eeas
Let $\dot{\Theta}$ be an open ball such that 
$\theta^*\in\dot{\Theta}$ and 
$\overline{\dot{\Theta}}\subset\Theta$. 
Since 
\beas
\sup_{\theta_0,\theta_1\in\dot{\Theta}}
\l|M_n(t;\theta_0,\theta_1)\r|\>|\theta_0-\theta_1|^{-1}
&\leq&
\sup_{\theta\in\dot{\Theta}} \Bigl| \calm_n(t;\theta) \Bigr|,
\eeas
one has
\beas 
P_{1,n}
&\leq&
P\bigg[\sup_{t:n\vartheta_n^2(t-t^*)>M}\frac{1}{[nt]-[nt^*]}
\l|M_n(t;\hat{\theta}_0,\hat{\theta}_1)\r|
\>|\hat{\theta}_0-\hat{\theta}_1|^{-1}
\geq\frac{c_\ep\vartheta_n}{6},\ 
%\\&& 
\hat{\theta}_0,\hat{\theta}_1\in\dot{\Theta}\bigg]
\\&&
+
P[|\hat{\theta}_0-\hat{\theta}_1|\geq 2\vartheta_n]
+
P[\hat{\theta}_0\not\in\dot{\Theta}]+P[\hat{\theta}_1\not\in\dot{\Theta}]
\\&\leq&
P\l[\sup_{t:n\vartheta_n^2(t-t^*)>M}\frac{1}{[nt]-[nt^*]}
\sup_{\theta\in\dot{\Theta}} \l|\calm_n(t;\theta)\r|
\geq\frac{c_\ep\vartheta_n}{6}\r]
\\&&
+
P[|\hat{\theta}_0-\hat{\theta}_1|\geq 2\vartheta_n]
+
P[\hat{\theta}_0\not\in\dot{\Theta}]+P[\hat{\theta}_1\not\in\dot{\Theta}]. 
\eeas

By the uniform version of 
the H\'ajek-Renyi inequality in Lemma \ref{unif.HR} 
applied to the case $p=2$, $s=2+[d_0/2]$ 
%where $d=1$, $p=2$, $s=1$ 
and $D=\dot{\Theta}$, 
we see under [H]$_2$ that 
\beas 
P\l[\sup_{t:n\vartheta_n^2(t-t^*)>M}\frac{1}{[nt]-[nt^*]}
\sup_{\theta\in\dot{\Theta}} \l|\calm_n(t;\theta)\r|
\geq\frac{c_\ep\vartheta_n}{6}\r]
&\leq&
\frac{{\sf C}}{c_\ep^2M}=:\rho_{\ep}(M), 
\eeas
therefore 
\bea\label{201221-11} 
\limsup_{n\to\infty}P_{1,n}
\leq 
\rho_{\ep}(M)
\eea
thanks to 
%\beas
%P[|\hat{\theta}_0-\hat{\theta}_1|\geq 2\vartheta_n]
%&\leq&
%P[|\hat{\theta}_0-\theta_0^*|\geq n^{-1/2}\ell_n]
%+
%P[|\hat{\theta}_1-\theta_1^*|\geq n^{-1/2}\ell_n]
%\\&&
%+
%1_{\{\vartheta_n<2n^{-1/2}\ell_n\}}.
%\eeas
\beas
P[|\hat{\theta}_0-\hat{\theta}_1|\geq 2\vartheta_n]
&\leq&
P[|\hat{\theta}_0-\theta_0^*|\geq \frac{1}{3}\vartheta_n]
+
P[|\hat{\theta}_1-\theta_1^*|\geq \frac{1}{3}\vartheta_n] 
\eeas
for large $n$.

In Case {\bf (A)}, 
Let $\dot{\Theta}_k$ be an open ball 
such that 
$\overline{\dot{\Theta}_k}\subset\Theta$ and 
$\theta_k^*\in
{\colord{\dot{\Theta}_k}}
%\ddot{\Theta}
$ 
for each $k=0,1$. 
\beas 
P_{1,n}&\leq&
P\l[\sup_{t:n\vartheta_0^2(t-t^*)>M}\frac{1}{[nt]-[nt^*]}
%\sup_{\theta_0,\theta_1\in\dot{\Theta}}
{\colord{
\sup_{
\theta_0\in\dot{\Theta}_0 \atop 
\theta_1\in\dot{\Theta}_1}
}}
\l|M_n(t;\theta_0,\theta_1)\r|
\geq\frac{c_\ep\vartheta_0^2}{3}\r]
\\&&
+
P[\hat{\theta}_0\not\in\dot{\Theta}_0]+P[\hat{\theta}_1\not\in\dot{\Theta}_1]. 
\eeas
We apply the H\'ajek-Renyi inequality for $M_n(t;\theta_0,\theta_1)$, 
which is a difference of two random fields on $\dot{\Theta}_k$ 
to be done with one by one, 
in order to obtain (\ref{201221-11}) under [H]$_1$.

\bigskip\noindent(iv) 
Estimation of $P_{2,n}$. 
First we consider Case {\bf (B)}. 
%There are positive constants $c_1$ and $c_2$ independent of $n$ 
There is a positive constant $c_2$ independent of $n$ 
such that 
\beas\label{201222-21} 
%\sum_{i=[nt^*]+1}^{[nt]} 
&&
{\rm Tr}\bigg(S_{i-1}(\hat{\theta}_0)^{-1}
S_{i-1}(\hat{\theta}_1)-I_d\bigg)
-\log\det\bigg(S_{i-1}(\hat{\theta}_0)^{-1}S_{i-1}(\hat{\theta}_1)\bigg) 
\nn\\&\geq&
\Xi(X_{t_{i-1}}^\tau,\theta^*)
[(\hat{\theta}_1-\hat{\theta}_0)^{\otimes2}]
+r_{n,i-1}|\hat{\theta}_1-\hat{\theta}_0|^2
\nn\\&\geq&
\{\lambda_1(\Xi(X_{t_{i-1}}^\tau,\theta^*))+r_{n,i-1}\}
|\hat{\theta}_1-\hat{\theta}_0|^2
\eeas
for all $i$, where $\max_i|r_{n,i-1}|\leq c_2\vartheta_n$, 
on the event 
\beas 
B_{K,n}&=&B_K\cap\{\hat{\theta}_0,\hat{\theta}_1\in\dot{\Theta},\>
|\hat{\theta}_k-\theta^*|\leq \vartheta_n\>(k=0,1)\}. 
\eeas
Thus 
\beas 
P_{2,n}
&\leq&
P\l[\inf_{t:n\vartheta_n^2(t-t^*)>M}\frac{1}{[nt]-[nt^*]}
A_n(t;\hat{\theta}_0,\hat{\theta}_1)\>|\hat{\theta}_1-\hat{\theta}_0|^{-2}
\leq 4c_\ep,\> B_{K,n}\r]
\\&&
+P\bigg[|\hat{\theta}_1-\hat{\theta}_0|\leq\half\vartheta_n\bigg]
+P[B_{K,n}^c]
\\&\leq&
P\l[\inf_{t:n\vartheta_n^2(t-t^*)>M}\frac{1}{[nt]-[nt^*]}
\sum_{i=[nt^*]+1}^{[nt]} 
\{\lambda_1(\Xi(X_{t_{i-1}}^\tau,\theta^*))+r_{n,i-1}\}
\leq 4c_\ep \r]
\\&&
+\ep
\eeas
for large $n$. 
%\mbox{\bf ここ検討}
%
\begin{en-text}
Here we used the identifiability inequality 
stemming from [A2-I] or [A2-II] on $B_K$. 
Due to the inequality, 
%\beas 
%P[|\hat{\theta}_0-\hat{\theta}_1|\leq \half\vartheta_n]
%&\leq&
%P[|\hat{\theta}_0-\theta_0^*|\geq n^{-1/2}\ell_n]
%+
%P[|\hat{\theta}_1-\theta_1^*|\geq n^{-1/2}\ell_n]
%\\&&
%+
%1_{\{\vartheta_n<4n^{-1/2}\ell_n\}},
%\eeas
\beas 
P[|\hat{\theta}_0-\hat{\theta}_1|\leq \half\vartheta_n]
&\leq&
P[|\hat{\theta}_0-\theta_0^*|\geq \frac{1}{6}\vartheta_n]
+
P[|\hat{\theta}_1-\theta_1^*|\geq \frac{1}{6}\vartheta_n]
\eeas
\end{en-text}
The scaled summation converges to the corresponding scaled integral 
uniformly in $t$ a.s., hence 
from Step (i) we have 
\beas 
\limsup_{n\iku\infty}P_{2,n}
%P\l[\inf_{t:n\vartheta_n^2(t-t^*)>M}\frac{1}{[nt]-[nt^*]}
%A_n(t;\theta_0^*,\theta_1^*)\leq c_\ep\vartheta_n^2\r]
%\\
&\leq&
P\l[
\inf_{t\in[t^*,1]}\frac{1}{t-t^*}\int_{t^*}^t 
%\l\{
%\frac{\sigma(s;\theta_1^*)^2}{\sigma(s;\theta_0^*)^2} - 1 -\log
%\frac{\sigma(s;\theta_1^*)^2}{\sigma(s;\theta_0^*)^2} \r\}
h(s)
\>ds \leq 
5c_\ep %\vartheta_n^2
\r]+\ep
\\&<&
2\ep
\eeas
for large $n$.

We will consider Case {\bf (A)}. 
There is a positive constant $c_2$ independent of $n$ 
such that 
\beas 
&&
{\rm Tr}\bigg(S_{i-1}(\hat{\theta}_0)^{-1}
S_{i-1}(\hat{\theta}_1)-I_d\bigg)
-\log\det\bigg(S_{i-1}(\hat{\theta}_0)^{-1}S_{i-1}(\hat{\theta}_1)\bigg) 
\\&\geq&
{\rm Tr}\bigg(S_{i-1}(\theta_0^*)^{-1}
S_{i-1}(\theta_1^*)-I_d\bigg)
-\log\det\bigg(S_{i-1}(\theta_0^*)^{-1}S_{i-1}(\theta_1^*)\bigg) 
\\&&
-
c_2(|\hat{\theta}_1-\theta_1^*|+|\hat{\theta}_0-\theta_0^*|)
\eeas 
for all $i$ on the event 
$
B'_{K,n}=B_K\cap\{\hat{\theta}_0\in\dot{\Theta}_0,
\hat{\theta}_1\in\dot{\Theta}_1\}
$ 
because there exists a continuous derivative $\partial_\theta\sigma$ 
by [H]$_1$. 
In this way, 
\beas 
P_{2,n}
&\leq&
P\l[\inf_{t:n(t-t^*)>M}\frac{1}{[nt]-[nt^*]}
A_n(t;\hat{\theta}_0,\hat{\theta}_1)
\leq c_\ep
{\colord{\vartheta_0^2}}
,\> B'_{K,n}\r]
\\&&
+P[B_{K,n}^{\prime c}]
\eeas
Therefore, 
\beas 
\limsup_{n\to\infty}P_{2,n} 
&\leq& 
P\l[
\inf_{t\in[t^*,1]}\frac{1}{t-t^*}\int_{t^*}^t 
%\l\{
%\frac{\sigma(s;\theta_1^*)^2}{\sigma(s;\theta_0^*)^2} - 1 -\log
%\frac{\sigma(s;\theta_1^*)^2}{\sigma(s;\theta_0^*)^2} \r\}
h(s)
\>ds \leq 
5c_\ep %\vartheta_n^2
\r]+\ep
\\&<&
2\ep
\eeas
by Step (i). 
\begin{en-text}
and the fact that 
$Q(X_{t^*},\theta_0^*,\theta_1^*)>0$ 
whenever $S(X_{t^*},\theta_0^*)\not=S(X_{t^*},\theta_1^*)$. 
\end{en-text}

\bigskip\noindent(v) 
Estimation of $P_{3,n}$. 
We have 
\beas
\sup_{t\in[t^*,1]} 
\l|
 S(X_t,\hat{\theta}_k)- S(X_t,\theta_k^*)\r|
\>1_{\{|\hat{\theta}_k-\theta_k^*|<2\vartheta_n\}\cap B_K}
\leq
{\sf C}\>\vartheta_n
\sskip(k=0,1),
\eeas
\beas
\sup_{t\in[t^*,1]} 
\l|
 S(X_t,\hat{\theta}_k)^{-1}- S(X_t,\theta_k^*)^{-1}
\r|
\>1_{\{|\hat{\theta}_k-\theta_k^*|<2\vartheta_n\}\cap B_K}
\leq
{\sf C}\>\vartheta_n
\sskip(k=0,1)
\eeas
and
\beas
\sup_{i:\geq [nt^*]+2} 
\l|S_{i-1}(\theta_1^*)
-h^{-1}E^{\theta_1^*}_{i-1}[(\Delta_iY)^{\otimes2}]\r|
1_{B_K}
&\leq&
{\sf C}\>
w_{[0,T]}(X,\frac{1}{n})^\alpha. 
\eeas
In the last estimate, the local $\alpha$-H\"older continuity 
of $\sigma$ was used. 
Then 
on $B_K\cap\{|\hat{\theta}_k-\theta_k^*|\leq2\vartheta_n\>(k=0,1)\}$, 
\bea\label{201223-10} 
&&
\sup_{t:n\vartheta_n^2(t-t^*)>M}\frac{1}{[nt]-[nt^*]}
\l|\rho_n(t;\hat{\theta}_0,\hat{\theta}_1)\r|
\vartheta_n^{-2}
=
o_p(1)
\eea
\begin{en-text}
\bea
\nn\\&\leq&
\sup_{t:n\vartheta_n^2(t-t^*)>M}\frac{{\sf C}}{[nt]-[nt^*]}
\sum_{i=[nt^*]+1}^{[nt]}
\bigg(|\hat{\theta}_1-\theta_1^*|
\nn\\&&
+h^{-1}\int_{t_{i-1}}^{t_i}[\sigma(X_{t_{i-1}}^\tau,\theta_1^*)^2-
\sigma(X_t^\tau,\theta_1^*)^2]\de t \bigg)
\vartheta_n^{-1}
\nn\\&\leq&
{\sf C}
\bigg(|\hat{\theta}_1-\theta_1^*|
\vartheta_n^{-1}
+o_p(1)\bigg) 
\eea
\end{en-text}
because of [H]$_j$ (ii). 
Consequently, we see 
$
\limsup_{n\to\infty}P_{3,n}\leq\ep$ due to [C] and 
the localization by $B_K$. 

\bigskip\noindent(vi) 
{From} the estimates in Steps (ii)-(iv) 
{\colorr and making $K$ sufficiently large}, 
we have 
\beas 
\limsup_{n\to\infty}
P[n{\colorr \vartheta_n^2}(\hat{t}_n-t^*)>M]
&\leq& 
\rho_{\ep}(M)+{\colorr 5}\ep
\eeas
for any $M\geq1$ and $\ep>0$. 
Therefore, 
\beas 
\limsup_{M\to\infty}\limsup_{n\to\infty}
P[n{\colorr \vartheta_n^2}(\hat{t}_n-t^*)>M]
&\leq& {\colorr 5}\ep, 
\eeas
which shows the tightness of 
$\{n{\colorr \vartheta_n^2}(\hat{t}_n-t^*)_+\}_n$. 
In a quite similar way, we can show that 
$\{n{\colorr \vartheta_n^2}(\hat{t}_n-t^*)_-\}_n$ 
is tight, and hence 
the family $\{n{\colorr \vartheta_n^2}(\hat{t}_n-t^*)\}_n$ 
is tight. \qed

%%%%%%%%%%%%%%%%%%%%%%%%%%%%%%%%%%%%%%%%%%%%%%%%%%%%%%%%%%%%
%%%%%%%%%%%%%%%%%%%%%%%%%%%%%%%%%%%%%%%%%%%%%%%%%%%%%%%%%%%%
%%%%%%%%%%%%%%%%%%%%%%%%%%%%%%%%%%%%%%%%%%%%%%%%%%%%%%%%%%%%

\section{Asymptotic distribution of 
{\colord{the change point estimator}}
}\label{sec:asym}
%\subsection{Case (B)}
%when $\theta_1^*$ converges to $\theta_0^*$} 

%Denote $\vartheta_n=\theta_1^*-\theta_0^*$. 

This section discusses limit theorems for the distributions of the estimators. 
First we consider Case {\bf (B)}.

Let 
\beas
\bbH(v)=
-2\left(
\Gamma_\eta^\half\>\calw(v)
-\half \Gamma_\eta|v|
\right)
\eeas
for $\Gamma_\eta={\colorr (}2{\colorr T)}^{-1}
\Xi(X_{t^*},\theta^*)[\eta^{\otimes2}]$. 
Here $\calw$ is a two-sided standard 
Wiener process independent of 
%$\sigma(X_{t^*},\theta^*)$. 
{\colord{$X_{t^*}$}}.

\begin{theorem}\label{201228-21}
%Let $T>0$. 
Suppose that the limit 
$\eta=\lim_{n\to\infty}\vartheta_n^{-1}(\theta_1^*-\theta_0^*)$ 
exists. 
Suppose that $[H]_2$, $[C]$ and $[B]$ are fulfilled in Case {\bf (B)}. 
Then 
$n\vartheta_n^2(\hat{t}_n-t^*)
%\to^{d_s(\calf_T)} 
{\colord{\to^d}}
{\rm argmin}_{v\in\bbR}\>\bbH(v)$ 
as $n\to\infty$. 
\end{theorem}

We will prove Theorem \ref{201228-21} 
and assume for a while that $T=1$ to simplify the notation. 
Introduce a new parameter $v$ as 
$t=t^\dagger_v:=t^*+v(n\vartheta_n^2)^{-1}$. 
Let 
\beas
D_n(v)
&=&
\l\{\Psi_n(t^\dagger_v;\hat{\theta}_0,\hat{\theta}_1)
-\Psi_n(t^*;\hat{\theta}_0,\hat{\theta}_1)\r\}
-
\l\{\Psi_n(t^\dagger_v;\theta_0^*,\theta_1^*)
-\Psi_n(t^*;\theta_0^*,\theta_1^*)\r\}
\\&=&
\{M_n(t^\dagger_v;\hat{\theta}_0,\hat{\theta}_1)
-M_n(t^\dagger_v;\theta_0^*,\theta_1^*)\}
+\{A_n(t^\dagger_v;\hat{\theta}_0,\hat{\theta}_1)
-A_n(t^\dagger_v;\theta_0^*,\theta_1^*)\}
\\&&
+\{\rho_n(t^\dagger_v;\hat{\theta}_0,\hat{\theta}_1)
-\rho_n(t^\dagger_v;\theta_0^*,\theta_1^*)\}. 
\eeas

\begin{lemma}\label{lemma-d}
For every $L>0$, 
\beas 
\sup_{v\in[-L,L]} |D_n(v)| \to^p 0
\eeas
as $n\to\infty$. 
\end{lemma}
\proof 
We assume that $v>0$. 
We have 
\beas 
&&
M_n(t^\dagger_v;\hat{\theta}_0,\hat{\theta}_1)
-M_n(t^\dagger_v;\theta_0^*,\theta_1^*)
\\&=&
\int_0^1 \vartheta_n \partial_\theta 
M_n(t^\dagger_v;\theta_0^*+u(\hat{\theta}_0-\theta_0^*), 
\theta_1^*+u(\hat{\theta}_1-\theta_1^*))\>du \>
[\vartheta_n^{-1}
(\hat{\theta}_0-\theta_0^*,\hat{\theta}_1-\theta_1^*)]. 
\eeas
For $k=0,1$ and $j=1,2$, 
\beas 
E\l[
\sup_{t\in [t^*,t^*+L(n\vartheta_n^2)^{-1}]}
|\partial_{\theta_k}^j M_n(t;\theta_0,\theta_1)|^2\r]
&\leq& 
{\colord{8}}E\l[
|\partial_{\theta_k}^j M_n(t^*+L(n\vartheta_n^2)^{-1};\theta_0,\theta_1)|^2\r]
{\colord{+O(1)}}
\\&\leq& 
{\colord{8}}L\vartheta_n^{-2}
\sup_{i\geq 1}
E\l[|\partial_{\theta_k}^j g_i(\theta_k)|^2\r]{\colord{+O(1)}}
\\&\leq&
{\sf C}L\vartheta_n^{-2}. 
\eeas
Then Sobolev's inequality implies 
%(now assume the dimension of $\theta_k$ is 
%one, consider $W^{1,2}\subset C$) 
\beas 
\vartheta_n \sup_{t\in [t^*,t^*+L(n\vartheta_n^2)^{-1}], 
\atop \theta_0,\theta_1\in\dot{\Theta}} 
|\partial_\theta M_n(t;\theta_0,\theta_1)|
=O_p(1).
\eeas
As a result, 
\begin{en-text}
\beas 
&&
\int_0^1 \vartheta_n \partial_\theta 
M_n(t^\dagger_v;\theta_0^*+u(\hat{\theta}_0-\theta_0^*), 
\theta_1^*+u(\hat{\theta}_1-\theta_1^*))\>du \>
[\vartheta_n^{-1}
(\hat{\theta}_0-\theta_0^*,\hat{\theta}_1-\theta_1^*)]
\\&=&
O_p(1)\times o_p(1)=o_p(1),
\eeas
and hence 
\end{en-text}
%This implies 
\beas 
\sup_{v\in[0,L]} 
|M_n(t^\dagger_v;\hat{\theta}_0,\hat{\theta}_1)
-M_n(t^\dagger_v;\theta_0^*,\theta_1^*)|
\to^p 0
\eeas
as $n\to\infty$. 

Set 
$r_n=|\hat{\theta}_0-\theta_0^*|+|\hat{\theta}_1-\theta_1^*|$. 
Simple calculus yields 
\beas 
|\{{\rm Tr}\>y-\log\det(I_d+y)\}-\{{\rm Tr}\>x-\log\det(I_d+x)\}|
&\leq& 
c_3 |y-x|(|x|+|y-x|)
\eeas
for $d\times d$-symmetric matrices $x$ and $y$ 
whenever $|x|,|y|\leq c_3'$, where $c_3'$ and 
$c_3$ are some positive constants independent of $x,y$. 
Indeed, the formula 
$\int\exp(-2^{-1}(I_d+\ep x)[z^{\otimes2}])dz
=(2\pi)^{d/2}\det(I_d+\ep x)^{-1/2}$ is convenient 
for explicite computation.

Applying this inequality to 
$y=S_{i-1}(\hat{\theta}_0)^{-1/2}S_{i-1}(\hat{\theta}_1)S_{i-1}(\hat{\theta}_0)^{-1/2}-I_d$ and 
$x=
%\sigma_{i-1}(\theta_0^*)^{-1/2}\sigma_{i-1}(\theta_1^*)\sigma_{i-1}(\theta_0^*)%^{-1/2}
{\colord{S_{i-1}(\theta_0^*)^{-1/2}S_{i-1}(\theta_1^*)S_{i-1}(\theta_0^*)^{-1/2}}}
-I_d$, 
we see that there exists a constant $c_4$ such that 
for large $n$, 
on $B_K\cap\{|\hat{\theta}_k-\theta^*|<\vartheta_n\>(k=0,1)\}$, 
\beas 
|A_n(t;\hat{\theta}_0,\hat{\theta}_1)-A_n(t;\theta_0^*,\theta_1^*)|
&\leq&
c_4\sum_{i=[nt^*]+1}^{[nt]}
r_n(\vartheta_n+r_n). 
\eeas
Thereofore, for any $\ep>0$, if we take sufficiently large $K$, then 
\beas 
\limsup_{n\to\infty}P\bigg[
\sup_{t\in [t^*,t^*+L(n\vartheta_n^2)^{-1}]}
|A_n(t;\hat{\theta}_0,\hat{\theta}_1)-A_n(t;\theta_0^*,\theta_1^*)|
\geq \ep \bigg]
\leq \ep. 
\eeas
This implies 
\beas 
\sup_{v\in[0,L]} 
|A_n(t^\dagger_v;\hat{\theta}_0,\hat{\theta}_1)
-A_n(t^\dagger_v;\theta_0^*,\theta_1^*)|
\to^p 0
\eeas
as $n\to\infty$. 
The convergence 
\beas 
\sup_{v\in[0,L]} 
|\rho_n(t^\dagger_v;\hat{\theta}_0,\hat{\theta}_1)
-\rho_n(t^\dagger_v;\theta_0^*,\theta_1^*)|
\to^p 0
\eeas
can be shown in the same way as (\ref{201223-10}). 

A similar proof of the uniform convergence on $[-L,0]$ 
is possible. 
After all, we obtained the desired result. \qed\\

\begin{remark}\rm 
When $\theta_k^*$ ($k=0,1$) are known, 
we do not need Lemma \ref{lemma-d}. 
\end{remark}

Thus we can focus only on $\Psi_n(t^\dagger_v;\theta_0^*,\theta_1^*)
-\Psi_n(t^*;\theta_0^*,\theta_1^*)$. 
For simplicity, we write $\Psi_n^*(t)$ for 
$\Psi_n(t;\theta_0^*,\theta_1^*)$. 
By assumption, %Assume that 
there exists a limit 
$\eta=\lim_{n\to\infty}\vartheta_n^{-1}(\theta_1^*-\theta_0^*)$. 
%
%The symbol $d_s(\calg)$ stands for the $\calg$-stable convergence 
%in distribution. 
$\bbD$ denotes the $D$-space on an interval of $t$. 
Let 
\beas 
\bbH_n(v)=\Psi_n^*\left(t^*+v(n\vartheta_n^2)^{-1}\right)-\Psi_n^*(t^*). 
\eeas
and 
\beas
\bbH^\tau(v)=
-2\left(
\Gamma_{\eta,\tau}^\half\>\calw(v)
-\half \Gamma_{\eta,\tau}|v|
\right)
\eeas
for $\Gamma_{\eta,\tau}={\colorr 1_{\{\tau>0\}}(}2{\colorr T)}^{-1}
\Xi({\colorr X_{t^*}^\tau},\theta^*)[\eta^{\otimes2}]$.

\begin{lemma}\label{201228-22}
Let 
$\eta=\lim_{n\to\infty}\vartheta_n^{-1}(\theta_1^*-\theta_0^*)$. 
Suppose that $[H]_2$, $[C]$ and $[B]$ are fulfilled in Case {\bf (B)}. 
Then 
$\bbH_n\to^{d_s(\calf_T)} \bbH^\tau$ in $\bbD([-L,L])$ as $n\to\infty$ 
for every $L>0$. 
\end{lemma}
\proof %{\bf Not yet completed. Corrections are necessary. }
\begin{en-text}
To demonstrate the procedure involving the negative values of 
the local parameter, 
we will consider $v\in[-L,0]$. 
\end{en-text}
We will only consider positive $v$ since 
the argument is essentially the same for negative $v$. 
Let $T=1$ as before. 
It follows from Lemma \ref{decomposition} that 
\beas 
\bbH_n(v)
&=& 
M_n^\Delta(v)
+A_n^\Delta(v)
+\rho_n^\Delta(v),
\eeas
where
\beas 
M_n^\Delta(v)
&=&
M_n(t^*+v(n\vartheta_n^2)^{-1};\theta_0^*,\theta_1^*), 
%-M_n(t^*;\theta_0^*,\theta_1^*),
\\
A_n^\Delta(v)
&=&
A_n(t^*+v(n\vartheta_n^2)^{-1};\theta_0^*,\theta_1^*),
%-A_n(t^*;\theta_0,\theta_1),
\\
\rho_n^\Delta(v)
&=&
\rho_n(t^*+v(n\vartheta_n^2)^{-1};\theta_0^*,\theta_1^*). 
%-\rho_n(t^*;\theta_0,\theta_1). 
\eeas
The evaluation of these terms will be done in the following. 
As repeated previously, we may proceed discussion on the event 
$B_K$ hereafter. 
First 
\bea\label{201228-10} 
M_n^\Delta(v)
&=&
1_{\{\tau>0\}}
\sum_{i=[nt^*]+1}^{[nt^*+\vartheta_n^{-2}v]} 
{\rm Tr} \bigg[ 
\big(S_{i-1}(\theta_0^*)^{-1}-S_{i-1}(\theta_1^*)^{-1}\big)
\nn\\&&\cdot 
h^{-1}\bigg(
\big(\int_{t_{i-1}}^{t_i}\sigma(X_t^\tau,\theta_1^*)dW_t\big)^{\otimes2} 
-E^{\theta_1^*}_{i-1} \big[
\int_{t_{i-1}}^{t_i}S(X_t^\tau,\theta_1^*)dt \big]
\bigg)\bigg]
\nn\\&&
+\bar{o}_p(1)
\eea
where $U_n(v)=\bar{o}_p(1)$ means that 
$\sup_{v\in[0,L]}|U_n(v)|\to^p0$, and we used the hypothesis 
$n\vartheta_n^2\to\infty$ and the fact that 
%$|\sigma_{i-1}(\theta_0^*)^{-2}-\sigma_{i-1}(\theta_1^*)^{-2}|
${\colord{
|S_{i-1}(\theta_0^*)^{-1}-S_{i-1}(\theta_1^*)^{-1}|
}}
\leq{\sf C}\>\vartheta_n$ with the localization. 
To obtain $\bar{o}_p(1)$, $L^1$-estimate helps. 
It follows from [H]$_j$ (i)(c) and (ii) that 
\beas
&&
\bigg|
h^{-1}\bigg(
\int_{t_{i-1}}^{t_i}S(X_t^\tau,\theta_1^*)dt
-E^{\theta_1^*}_{i-1}\bigg[
\int_{t_{i-1}}^{t_i}S(X_t^\tau,\theta_1^*)dt\bigg]
\bigg)\bigg|
\\&=&
\bigg|
h^{-1}\bigg(
\int_{t_{i-1}}^{t_i}
[S(X_t^\tau,\theta_1^*)-S(X_{t_{i-1}}^\tau,\theta_1^*)]dt
\\&&
-E^{\theta_1^*}_{i-1}\bigg[
\int_{t_{i-1}}^{t_i}
[S(X_t^\tau,\theta_1^*)-S(X_{t_{i-1}}^\tau,\theta_1^*)]dt\bigg]
\bigg)\bigg|
\\&\leq&
{\sf C}\> w_{[0,T]}(n^{-1},X)^\alpha
\\&=&
{\colord{o_p(\vartheta_n)}}.
\eeas
Moreover, with the Burkholder-Davis-Gundy inequality, 
\begin{en-text}
we obtain
\beas 
&&1_{\{\tau>0\}}
\sum_{i=[nt^*]+1}^{[nt^*+\vartheta_n^{-2}v]} 
{\rm Tr}\bigg[
(\sigma_{i-1}(\theta_0^*)^{-2}-\sigma_{i-1}(\theta_1^*)^{-2})
\\&&\cdot 
h^{-1}\bigg(
(\int_{t_{i-1}}^{t_i}\sigma(X_t^\tau,\theta_1^*)dW_t)^2 
-
\int_{t_{i-1}}^{t_i}\sigma(X_t^\tau,\theta_1^*)^2dt
\bigg)
\\&=&
\tilde{M}_n^\Delta(v)+\bar{o}_p(1), 
\eeas
\end{en-text}
{\colord{the first terms on}} 
the right-hand side of (\ref{201228-10}) equals 
{\colord{
$\bar{M}_n^\Delta(v)+\bar{o}_p(1)$ with 
$\bar{M}_n^\Delta(v)=
\sum_{i=[nt^*]+1}^{[nt^*+\vartheta_n^{-2}v]} \xi_{n,i}$}}, 
where
\bea\label{210124-12} 
\xi_{n,i}
&=&
1_{\{\tau>0\}}
%\sum_{i=[nt^*]+1}^{[nt^*+\vartheta_n^{-2}v]} 
{\rm Tr}\bigg[ 
{^{\sf t}}\sigma_{i-1}(\theta_1^*)
\big(S_{i-1}(\theta_0^*)^{-1}-S_{i-1}(\theta_1^*)^{-1}\big)
\sigma_{i-1}(\theta_1^*)
\nn\\&&\cdot 
\big(h^{-1}(\Delta_iW)^{\otimes2} -I_r\big)
\bigg]
\eea 
{\colord{and 
$\sigma_{i-1}(\theta)=\sigma(X_{t_{i-1}}^\tau,\theta)$. 

Now we introduce the backward approximation 
\beas 
\tilde{\xi}_{n,i}
&=&
1_{\{\tau>0\}}
{\rm Tr}\bigg[ 
{^{\sf t}}\sigma(X_{t^*-\ep_n}^\tau,\theta_1^*)
\big(S(X_{t^*-\ep_n}^\tau,\theta_0^*)^{-1}
-S(X_{t^*-\ep_n}^\tau,\theta_1^*)^{-1}\big)
\sigma(X_{t^*-\ep_n}^\tau,\theta_1^*)
\\&&\cdot 
\big(h^{-1}(\Delta_iW)^{\otimes2} -I_r\big)
\bigg]
\eeas
to $\xi_{n,i}$ 
for $\ep_n=2Ln^{-1}\vartheta_n^{-2}$. 
}}
After all, 
{\colord{by 
$\tilde{M}_n^\Delta(v)=
\sum_{i=[nt^*]+1}^{[nt^*+\vartheta_n^{-2}v]} \tilde{\xi}_{n,i}$, 
we have}}
\bea\label{210124-11} 
M_n^\Delta(v)
&=&
\tilde{M}_n^\Delta(v)+\bar{o}_p(1), 
\eea
\begin{en-text}
Let $\calg^n_v=\calf_{[nt^*+\vartheta_n^{-2}v]h}$. 
The predictable quadratic process of $\tilde{M}_n^\Delta$ 
with respect to the filtration 
$\G^n=(\calg^n_v)$ 
is given by 
\beas 
&&
\langle\tilde{M}_n^\Delta\rangle(v)
\\&=&
1_{\{\tau>0\}}
\sum_{i=[nt^*]+1}^{[nt^*+\vartheta_n^{-2}v]} 2
%
%\big|{^{\sf t}}\sigma_{i-1}(\theta_1^*)
%(S_{i-1}(\theta_0^*)^{-1}-S_{i-1}(\theta_1^*)^{-1})
%\sigma_{i-1}(\theta_1^*)\big|^2
\big|
{\colord{
{^{\sf t}}\sigma(X_{t^*-\ep_n}^\tau,\theta_1^*)
\big(S(X_{t^*-\ep_n}^\tau,\theta_0^*)^{-1}
-S(X_{t^*-\ep_n}^\tau,\theta_1^*)^{-1}\big)
\sigma(X_{t^*-\ep_n}^\tau,\theta_1^*)
}}
\big|^2
%
%(\sigma_{i-1}(\theta_0^*)^{-2}-\sigma_{i-1}(\theta_1^*)^{-2})^2
%\sigma_{i-1}(\theta_1^*)^4
%(h^{-1}(\Delta_iW)^2 -1)^2
\\&=&
1_{\{\tau>0\}}
\sum_{i=[nt^*]+1}^{[nt^*+\vartheta_n^{-2}v]} 2
%\vartheta_n^2
%(\partial_\theta\sigma_{i-1}(\theta^*)/\sigma_{i-1}(\theta^*))^2
\Xi(X_{{\colord{t^*-\ep_n}}}^\tau,\theta^*)[(\theta_1^*-\theta_0^*)^{\otimes2}]
+\bar{o}_p(1)
\\&\to^p&
1_{\{\tau>0\}}
2
%(\partial_\theta\sigma/\sigma(X_{t^*},\theta^*))^2v.
\Xi(X^\tau_{t^*},\theta^*)[\eta^{\otimes2}]\>v. 
\eeas
\end{en-text}
\begin{en-text}
Besides, 
\beas 
\langle \tilde{M}_n^\Delta(v), W_n \rangle(v)
&=&
0
\eeas
for $W_n(v)=W_{[nt^*+\vartheta_n^{-2}v]h}$. 
For any bounded martingale $N$ orthogonal to $W$, 
if we set $N_n(v)=N_{[nt^*+\vartheta_n^{-2}v]h}$, then 
\beas 
\langle \tilde{M}_n^\Delta(v), N_n \rangle(v)
&=&
0.
\eeas
\end{en-text}
\begin{en-text}
Also, 
it is easy to see that for any $\ep>0$, 
\beas 
\sum_{i=[nt^*]+1}^{[nt^*+\vartheta_n^{-2}v]} 
E_{i-1}^{\theta_1^*}[
{\colord{\tilde{
\xi}_{n,i}^2
}}
1_{\{|\xi_{n,i}|>\ep\}}]
&\to^p&
0. 
\eeas
\end{en-text}
\begin{en-text}
for 
\beas 
\xi_{n,i}
&=&
(\sigma_{i-1}(\theta_0^*)^{-2}-\sigma_{i-1}(\theta_1^*)^{-2})
\sigma_{i-1}(\theta_1^*)^2
%\\&&
(h^{-1}(\Delta_iW)^2 -1). 
\eeas 
\end{en-text}
%
%
%Then Theorem 3-2 of Jacod (1997) 
{\colord{
Since 
\beas 
&&
1_{\{\tau>0\}}
%\sum_{i=[nt^*]+1}^{[nt^*+\vartheta_n^{-2}v]} 
2\vartheta_n^{-2}v
%
%\big|{^{\sf t}}\sigma_{i-1}(\theta_1^*)
%(S_{i-1}(\theta_0^*)^{-1}-S_{i-1}(\theta_1^*)^{-1})
%\sigma_{i-1}(\theta_1^*)\big|^2
\big|
{\colord{
{^{\sf t}}\sigma(X_{t^*-\ep_n}^\tau,\theta_1^*)
\big(S(X_{t^*-\ep_n}^\tau,\theta_0^*)^{-1}
-S(X_{t^*-\ep_n}^\tau,\theta_1^*)^{-1}\big)
\sigma(X_{t^*-\ep_n}^\tau,\theta_1^*)
}}
\big|^2
%
%(\sigma_{i-1}(\theta_0^*)^{-2}-\sigma_{i-1}(\theta_1^*)^{-2})^2
%\sigma_{i-1}(\theta_1^*)^4
%(h^{-1}(\Delta_iW)^2 -1)^2
\\&=&
1_{\{\tau>0\}}
%\sum_{i=[nt^*]+1}^{[nt^*+\vartheta_n^{-2}v]} 
2\vartheta_n^{-2}v 
%\vartheta_n^2
%(\partial_\theta\sigma_{i-1}(\theta^*)/\sigma_{i-1}(\theta^*))^2
\Xi(X_{{\colord{t^*-\ep_n}}}^\tau,\theta^*)[(\theta_1^*-\theta_0^*)^{\otimes2}]
+\bar{o}_p(1)
\\&\to^p&
1_{\{\tau>0\}}
2
%(\partial_\theta\sigma/\sigma(X_{t^*},\theta^*))^2v.
\Xi(X^\tau_{t^*},\theta^*)[\eta^{\otimes2}]\>v, 
\eeas
\begin{en-text}
the martingale central limit theorem of mixture type 
or the theory of the convergence of stochastic integrals 
\end{en-text}
the central limit theorem}} 
ensures {\colord{the convergence}} %the stable convergence 
%$\tilde{M}_n^\Delta\to^{d_s(\calf_1)}-2\Gamma_{\eta,\tau}^\half\>\calw$
{\colord{
$\tilde{M}_n^\Delta\to^{d}
-2\Gamma_{\eta,\tau}^\half\>\calw$
}}
in $\bbD([0,L])$. 
{\colord{
Indeed, the joint convergence 
of 
$X^\tau_{t^*-\ep_n}$, 
$\vartheta_n^{-1}
{^{\sf t}}\sigma(X_{t^*-\ep_n}^\tau,\theta_1^*)
\big(S(X_{t^*-\ep_n}^\tau,\theta_0^*)^{-1}
-S(X_{t^*-\ep_n}^\tau,\theta_1^*)^{-1}\big)
\sigma(X_{t^*-\ep_n}^\tau,\theta_1^*)$ 
and the process 
$\vartheta_n\sum_{i=[nt^*]+1}^{[nt^*+\vartheta_n^{-2}v]} 
\big(h^{-1}(\Delta_iW)^{\otimes2} -I_r\big)$ 
implies the joint convergence of 
$(X^\tau_{t^*},\tilde{M}^\Delta_n)$. 
}}
In the same fashion, we can show 
$\tilde{M}_n^\Delta\to^d %{d_s(\calf_1)}
-2\Gamma_{\eta,\tau}^\half\>\calw'(\cdot+L)$ in $\bbD([-L,0])$ 
if $\tilde{M}_n^\Delta$ is defined in a natural way 
over negative $v$, where 
$(\calw'(u))_{u\in[0,L]}$ is a standard Wiener process 
independent of $(\calw(v))_{v\in[0,L]}$ and $\calf$. 
Since $\Xi(X_{t^*}^\tau,\theta^*)[\eta^{\otimes2}]$ is 
independent of $\calw'$, we can replace 
the stochastic integral with respect to $\calw'$ in 
the representation of the limit distribution of 
$\tilde{M}_n^\Delta$ by the one with respect to 
the negative-time part of the two sided Wiener process 
$\calw$ reversible in time. 
Easy calculations yield
$\sup_{v\in[-L,L]}|A_n^\Delta(v)-\Gamma_{\eta,\tau} v|\to^p0$ and
$\sup_{v\in[-L,L]}|\rho_n^\Delta(v)|\to^p0$ 
for extended $A_n^\Delta$ and $\rho_n^\Delta$ to $[-L,L]$, 
which completes the proof. \qed \\

{\it Proof of Theorem \ref{201228-21}.} 
We have supposed that $T=1$ to state the lemmas, and we start with this case. 
Write $\hat{v}={\rm argmin}_{v\in\bbR}\>\bbH(v)$. 
For $\ep>0$, take large $K$ so that $P[\tau=T]>1-\ep$. 
It follows from 
Lemma \ref{201228-22} that 
for every $x\in\bbR$, 
\beas &&
\limsup_{n\to\infty}
P[n\vartheta_n^2(\hat{t}_n-t^*)\leq x]-\ep
\\&\leq& 
\limsup_{n\to\infty}
P[\inf_{v\in[-L,x]}\bbH_n^\tau(v)\leq \inf_{v\in[x,L]}\bbH_n^\tau(v)]
%\\&&
+\sup_nP[n\vartheta_n^2(\hat{t}_n-t^*)\not\in[-L,L]\}
\\&=&
P[\inf_{v\in[-L,x]}\bbH^\tau(v)\leq \inf_{v\in[x,L]}\bbH^\tau(v)]
%\\&&
+\sup_nP[n\vartheta_n^2(\hat{t}_n-t^*)\not\in[-L,L]\}
\\&\leq&
\ep+P[\hat{v}\leq x]
%\\&&
+P[\hat{v}\not\in[-L,L]]
+\sup_nP[n\vartheta_n^2(\hat{t}_n-t^*)\not\in[-L,L]] 
\eeas 
As $L\to\infty$, the last two terms of the right-hand side 
of the above inequality tend to $0$ 
thanks to Theorem \ref{081219-1} (b). 
So we have obtained 
\beas 
\limsup_{n\to\infty}
P[n\vartheta_n^2(\hat{t}_n-t^*)\leq x]
&\leq&
P[\hat{v}\leq x]. 
\eeas
The estimate of $P[n\vartheta_n^2(\hat{t}_n-t^*)\leq x]$ 
from below can be done in a similar manner, which 
concludes the proof in case $T=1$. 

For general $T$, 
we introduece 
a stochastic basis 
$\tilde{\calb}=(\Omega,\calf,\tilde{\F},P)$ with 
$\tilde{\F}=(\calf_{Tu})_{u\in[0,1]}$, and 
the processes 
$\tilde{b}_u=b_{Tu}$, $\tilde{X}_u=X_{Tu}$ and 
$\tilde{Y}_u=Y_{Tu}$, $u\in[0,1]$, to scale the time as $t=Tu$. 
Those stochastic processes satisfy the stochastic integral equation 
\beas 
\tilde{Y}_u &=& \tilde{Y}_0 + 
\int_0^u \tilde{b}_rdr+\int_0^u \tilde{\sigma}(\tilde{X}_r,\theta)
d\tilde{W}_r,
\eeas
where $\tilde{\sigma}(x,\theta)=\sqrt{T}\sigma(x,\theta)$ and 
$\tilde{W}$ is an $r$-dimensional $\tilde{\F}$-Wiener process. 
The sampling times $(iT/n)_{i=0}^n$ now chage to $(i/n)_{i=0}^n$ 
in the new setting after scaling time. 
For the change point estimator $\hat{u}_n$ for $u^*=T^{-1}t^*$, 
we know 
\bea\label{210123-1} 
n\vartheta_n^2(\hat{u}_n-u^*)&\to^{d_s}& 
\mbox{argmin}_{\tilde{v}\in\bbR}\> \tilde{\bbH}(\tilde{v}),
\eea
where 
$
\tilde{\bbH}(\tilde{v}) = 
-2\big(\tilde{\Gamma}_\eta\tilde{\calw}(\tilde{v})
-2^{-1}\tilde{\Gamma}_\eta|\tilde{v}|\big) 
$, 
$\tilde{\Gamma}_\eta=2^{-1}\Xi(\tilde{X}_{u^*},\theta^*)[\eta^{\otimes2}]$ and 
$\tilde{\calw}$ is a two-sided Wiener process independent of 
$\tilde{\sigma}(\tilde{X}_{u^*},\theta^*)=\sqrt{T}\sigma(X_{t^*},\theta^*)$. 
Since 
\beas 
T\>\mbox{argmin}_{\tilde{v}\in\bbR}\> \tilde{\bbH}(\tilde{v})
&=&
\mbox{argmin}_{v\in\bbR}\> \tilde{\bbH}\left(\frac{v}{T}\right)
\\&=^d&
\mbox{argmin}_{v\in\bbR}\>\bbH(v)
\eeas
thanks to 
$\calw(\cdot)=^d T^{1/2}\tilde{\calw}(\cdot/T)$. 
Thus (\ref{210123-1}) gives the desired convergence of $\hat{t}_n$ 
since $\hat{t}_n=T\hat{u}_n$. 
\qed \\

Let us investigate the limit distribution of 
the estimator in Case {\bf (A)}. 
By nature of the sampling scheme, only 
the set 
$\calg_n=\{kT/n;k\in\bbZ\}$ has essential meaning 
for the optimization with respect to 
the parameter $t$. Without loss of generality, 
we modify $\hat{t}_n$ so that it takes values in $\calg_n$, and 
set $\hat{k}_n=n\hat{t}_n/T$. 
%In the same sense, the true value $t^*$ only has 
Let
\beas 
\bbK(v)= 
&&
%\sum_{i=[nt^*/T]+1}^{[nt^*/T]+v} \bigg\{
\sum_{i=1}^{v} \bigg\{
{\rm Tr}\bigg[ 
{^{\sf t}}\sigma(X_{t^*},\theta_1^*)
\big(S(X_{t^*},\theta_0^*)^{-1}
-S(X_{t^*},\theta_1^*)^{-1}\big)
\sigma(X_{t^*},\theta_1^*)
%\\&&\cdot \big(
\zeta_i^{\otimes2} %-I_r\big)
\bigg]
%+
%{\rm Tr}\bigg(S(X_{t^*},\theta_0)^{-1}
%S(X_{t^*}^\tau,\theta_1)-I_d\bigg)
\\&& 
-\log\det\bigg(S(X_{t^*},\theta_0^*)^{-1}
S(X_{t^*},\theta_1^*)\bigg) 
\bigg\},
\eeas
where $\zeta_i$ are independent $r$-dimensional standard normal 
variables independent of $X_{t^*}$. 
\begin{theorem}\label{210124-20}
Suppose that $[H]_1$, $[C]$ and $[A]$ are fulfilled in Case {\bf (A)}. 
Then 
$\hat{k}_n-[\frac{nt^*}{T}]
%=\frac{n}{T}(\hat{t}_n-[\frac{nt^*}{T}]\frac{T}{n})
\to^d 
{\rm argmin}_{v\in\bbZ}\>\bbK(v)$ 
as $n\to\infty$. 
\end{theorem}
\proof 
We change the definition of $t^\dagger$ and 
newly set 
$t^\dagger_v=[\frac{nt^*}{T}]\frac{T}{n}+\frac{Tv}{n}$. 
Lemma \ref{lemma-d} is still valid 
by essentially the same proof and hence we may only consider 
$\Psi_n(t^\dagger_v;\theta_0^*,\theta_1^*)
-\Psi_n(t^*;\theta_0^*,\theta_1^*)$. 
Writing $\Psi_n^*(t)$ for 
$\Psi_n(t;\theta_0^*,\theta_1^*)$, 
we will investigate the behavior of the random field 
\beas 
\bbK_n(v)
&=&
\Psi_n^*(t^\dagger_v)-\Psi_n^*(t^*)
\eeas
on $v\in\bbZ$. 
For a while, 
we consider nonnegative $v$. The argument is similar for negative $v$. 
According to Lemma \ref{decomposition}, 
we have the decomposition 
\beas 
\bbK_n(v) = \bbM_n(v)+\bbA_n(v)+\varrho_n(v),
\eeas
where 
$
\bbM_n(v) = M_n(t^\dagger_v;\theta_0^*,\theta_1^*)$, 
$\bbA_n(v) = A_n(t^\dagger_v;\theta_0^*,\theta_1^*)$ and 
$\varrho_n(v) = \rho_n(t^\dagger_v;\theta_0^*,\theta_1^*)$.  

Now, $\bbM_n(v)$ admits a similar expansion 
as {\colord{before}}:%(\ref{210124-11}):
\beas 
\bbM_n(v)
&=&
\sum_{i=[nt^*/T]+1}^{[nt^*/T]+v} \xi_{n,i}+\bar{o}_p(1)
\eeas
with $\xi_{n,i}$ given by (\ref{210124-12}). 
Moreover, for $\ep_n=n^{-1/2}$ this time, we consider 
the backward approximation of $\xi_{n,i}$, 
{\colord{that is, }}
\beas 
\xi_{n,i}
&=&
\tilde{\xi}_{n,i}+\bar{o}_p(1). 
\eeas
Here $v\in[0,L]\cap\bbZ$, however this approximation 
is available when we consider $v\in[-L,0]$. 
let $L_0$ be the maximum integer in $[0,L]$. 
By continuity of $\sigma$ and because $W$ is an $\F$-Wiener process, 
we have 
\beas 
%\big(X_{t^*},(h^{-1}(\Delta_iW)^{\otimes2})
%_{i=[nt^*/T]-L_0}^{[nt^*/T]+L_0}\big)
%&=& 
\big(X_{t^*-\ep_n}^\tau,(h^{-1}(\Delta_iW)^{\otimes2})
_{i=[nt^*/T]-L_0}^{[nt^*/T]+L_0}\big)
%+o_p(1)
&\to^d& 
\big(X_{t^*}^\tau,(\zeta_i^{\otimes2})_{i=-L_0}^{L_0}\big),
\eeas
where $\zeta_i$ are independent $r$-dimensional standard normal 
variables independent of $X_{t^*}^\tau$; we use the same symbol $\zeta_i$ 
as in the statement. 
Consequently, 
\beas 
(X_{t^*}^\tau,\bbM_n(v))_{v=-L_0}^{L_0}
\to^d
(X_{t^*}^\tau,\bbM_\infty(v))_{v=-L_0}^{L_0},
\eeas 
where 
\beas 
\bbM_\infty(v)
&=&
\sum_{i=[nt^*/T]+1}^{[nt^*/T]+v} \xi_{\infty,i}
\eeas
and $\xi_{\infty,i}$ is given by 
\beas 
\xi_{\infty,i}
&=&
1_{\{\tau>0\}}
{\rm Tr}\bigg[ 
{^{\sf t}}\sigma(X_{t^*}^\tau,\theta_1^*)
\big(S(X_{t^*}^\tau,\theta_0^*)^{-1}
-S(X_{t^*}^\tau,\theta_1^*)^{-1}\big)
\sigma(X_{t^*}^\tau,\theta_1^*)
\\&&\cdot 
\big(\zeta_i^{\otimes2} -I_r\big)
\bigg]. 
\eeas

For $\bbA_n$, we have 
$\bbA_n(v)\to\bbA_\infty(v)$ with 
\beas
\bbA_\infty(v)
&=& 
1_{\{\tau>0\}}
\sum_{i=[nt^*/T]+1}^{[nt^*/T]+v} 
\bigg\{
{\rm Tr}\bigg(S(X_{t^*}^\tau,\theta_0^*)^{-1}
S(X_{t^*}^\tau,\theta_1^*)-I_d\bigg)
\\&& 
-\log\det\bigg(S(X_{t^*}^\tau,\theta_0^*)^{-1}
S(X_{t^*}^\tau,\theta_1^*)\bigg) 
\bigg\}.
\eeas

On the other hand, $\varrho_n(v)$ tends to $0$ uniformly in $v$. 
Therefore, 
\beas 
\big(\bbK_n(v)\big)_{v=-L_0}^{L_0} 
&\to^d&
\big(\bbK^\tau(v)\big)_{v=-L_0}^{L_0}, 
\eeas  
where $\bbK^\tau(v)=\bbM_\infty(v)+\bbA_\infty(v)$. 
Removing $\tau$ by letting $K\to\infty$, and 
using Theorem \ref{081219-1}, we obtain the limit distribution of $\hat{t}_n$. 
\qed\\

\begin{en-text}
and after approxiamtion as in the proof of
Lemma \ref{201228-22}, the problem reduces in showing
\beas 
\bigg(
h^{-1}
\big(\int_{t_{i-1}}^{t_i} dW_t \big)^{\otimes2}-I_r
\bigg)_{i=0,...,m}
\to^{d_s(\calf)}
\big(\zeta_i^{\otimes2}-I_r\big)_{i=0,...,m}
\eeas 
for independent $r$-dimensional standard normal random variables 
$\zeta_0,...,\zeta_m$ that are independent of $\calf$ 
or at least $X_{t^*}$, for every $m$. 
{\bf [??????????????????????????]}
At least we need independency of $\zeta_i$'s and $X_{t^*}$. 
When $v$ is positive, the independency is clear. 
Even when $v$ is nonpositive, the number of $v$'s is finite, 
so we can first replace $X_{t^*}$ by $X_{t^*-N/n}$ and 
show independency of $\zeta_i$'s and $\calf_{t^*-\ep}$, and 
finally take limit for changing $X_{t^*-\ep}$ by $X_{t^*}$ 
in the expression of the limit, 
by using the continuity of the expression. 
\end{en-text}

\section{Initial estimator for $\theta_k$}\label{081227-1} 

In this section, we will breifly discuss 
the consturction of the initial estimators. 
There are two situations according to the prior knowledge 
of the parameter space $\bbT$ of the change point. 
The first one is the case where 
$\bbT=[{\sf t}_0,{\sf t}_1]\subset(0,1)$ for given numbers 
${\sf t}_0$ and ${\sf t}_1$. 
In the second case, we do not assume a prior information 
of ${\sf t}_0$ and ${\sf t}_1$, instead the precision of the initial 
estimator will be lost. 
Let 
\beas 
\Phi_n^{0}(t;\theta_0)
=
\sum_{i=1}^{[nt{\colorr /T}]}G_i(\theta_0) 
\sskip\mbox{and}\sskip 
\Phi_n^{1}(t;\theta_1)
=
\sum_{i={[nt{\colorr /T}]+1}}^nG_i(\theta_1). 
\eeas
Suppose that ${\sf t}_0$ and ${\sf t}_1$ are known. 
Let $\hat{\theta}_0$ and $\hat{\theta}_1$ 
satisfy 
\beas 
\Phi_n^{k}({\sf t}_k;\hat{\theta}_k)
=
\min_{\theta_k}\Phi_n^{k}({\sf t}_k;\theta_k) 
\eeas
for $k=0,1$. 
To validate asymptotic properties of the estimators, 
it is sufficient that these relations are satisfied 
asymptotically. 
Under suitable regularity conditions as well as 
the identifiability conditions that 
\bea\label{210312-3} 
\int_0^{{\sf t}_0}
Q(X_t,\theta^*,\theta)
\>dt
>0\sskip a.s. 
%\eea
&\mbox{ and }&
%\bea 
\int_{{\sf t}_1}^T
Q(X_t,\theta^*,\theta)
\>dt
>0\sskip a.s. 
\eea
for every $\theta\not=\theta^*$, 
it is possible to show that 
$\hat{\theta}_k-\theta_k^*=O_p(n^{-1/2})$, 
therefore 
Condition [C] is satisfied in both cases {\bf (A)} and {\bf (B)}. 
Based on $\hat{\theta}_k$, the estimator $\hat{t}_n$ 
are defined. 
According to the previous sections, $\hat{t}_n$ 
possesses $n\vartheta_n^2$-consistency and 
the asymptotic distribution in each case is already known. 

We can also construct the second stage estimators. 
Let $b_n$ be a sequence of positive numbers such that 
{\colorr{$b_n=(n\vartheta_n^\delta)^{-1}$, 
where $\delta\in(2,\infty)$ is a constant satisfying 
$n\vartheta_n^\delta\to\infty$ as $n\to\infty$. }}
Construct $\check{\theta}_k$ so that 
\beas 
\Phi_n^{k}(\hat{t}_n+(-1)^{k+1}b_n;\check{\theta}_k)
=
\min_{\theta_k}\Phi_n^{k}(\hat{t}_n+(-1)^{k+1}b_n;\theta_k) 
\eeas
for $k=0,1$. 
The new estimators $\check{\theta}_k$ are expected 
to improve $\hat{\theta}_k$ since they utilize up to 
the data near $t^*$. 
{\colord{Further, it is possible to construct 
a new change-point estimator with those estimators. }}
Based on $\check{\theta}_k$, we define $\check{t}_n$ for $t^*$ 
as 
\beas 
\check{t}_n&=&\mbox{argmin}_{t\in[0,T]} 
\Phi_n(t;\check{\theta}_0,\check{\theta}_1). 
\eeas
Since 
it is usually easy to verify 
Condition [C] for $\check{\theta}_k$, 
we will then obtain the same asympotic results for $\check{t}_n$ 
as $\hat{t}_n$. 
%{\bf Write} 

\vspace*{3mm} 
Next, let us consider the second situation. 
The knowledge of ${\sf t}_k$ is not available 
and it means that any data set sampled over a fixed time interval 
$[0,a]$ is usuless for estimating $\theta_0$ 
since $t^*$ may be less than $a$ and then 
the data over $(t^*,a]$ causes bias in general. 
A similar notice is also for the estimation of $\theta_1$. 
This consideration suggests the use of estimators $\hat{\theta}_k$ 
based on the data over time interval $[0,a_n]$ for $k=0$ and 
the one over $[T-a_n,T]$ for $k=1$, respectively, 
for some sequence $a_n$ tending to zero. 
We assume that 
there exitst a constant $\beta\in(0,1/2)$ such that 
$a_n\geq 1/(n\vartheta_n^{1/\beta})$ and that 
\begin{en-text}
\beas 
\lim_{n\to\infty}
P\big[|\hat{\theta}_k-\theta_k^*|\geq (na_n)^{-\beta}\big]=0
\sskip(k=0,1). 
\eeas
\end{en-text}
$|\hat{\theta}_k-\theta_k^*|=o_p((na_n)^{-\beta})$ for $k=0,1$. 
When $\limsup_{n\to\infty}\vartheta_n>0$, we also assume $na_n\to\infty$. 
In particular, the first condition implies $n\vartheta_n^2\to\infty$. 
The second condition is natural because the number of data 
is proportional to $na_n$. 
To obtain $\hat{\theta}_k$, we may need the identifiability condition that 
$\sigma(\theta,x)=\sigma(\theta',x)$ implies $\theta=\theta'$; 
it is a strong condition 
like monotonicity of $\sigma(\theta,x)$ in $\theta$. 
Under the assumptions, [C] holds and after that 
it is possible to construct $\hat{t}_n$, $\check{\theta}_k$ and 
$\check{t}_n$ in turn as mentioned above. 
The asymptotic properties of $\check{t}_n$ 
are the same as $\hat{t}_n$ 
because $\check{\theta}_k$'s satisfy Condition [C]. 
It is expected that the new estimator $\check{t}_n$ posseses 
equal or better precision than $\hat{t}_n$ 
as numerical studies in Section \ref{201229-1} suggest.

%%%%%%%%%%%%%%%%%%%%%%%%%%%%%%%%%%%%%%%%%%%%%%%%%%%%%%%%%%%%%%%%%%%%%
\begin{en-text}
This section gives an example of 
initial estiamtor estimator for $\theta_k$. 
$\hat{\theta}_0([0,a_n])$ is 
the estimator constructed by the data 
discretely sampled from the interval $[0,a_n]$. 
$a_n\to0$ as $n\to\infty$. 
Let $\hat{\theta}_0=\hat{\theta}_0([0,a_n])$. 

It is possible to improve $\hat{\theta}_k$ if we 
take advantage of $\hat{t}^*$ based on 
$\hat{\theta}_k$. 
Define $\check{\theta}_k$ by 
\beas 
\check{\theta}_0 = \hat{\theta}_0([0,\hat{t}_n-a'_n]) 
\eeas
and 
{\bf koko}

\begin{theorem}
Then $\check{\theta}_k$ is 
a $\sqrt{n}$-consistent estimator for $\theta_k$. 
\end{theorem}

The second-stage estimator for $t^*$ is 
given by $\check{t}_n$ 
as the minimum point of $\Phi_n(t;\check{\theta}_0,\check{\theta}_1)$. 

\begin{theorem}
\begin{description}
\item[(a)] 
asymptotic distribution of $\check{\theta}_n$ in Case {\bf (A)}
\item[(b)] asymptotic mixed normality of $\check{\theta}_n$ 
in Case {\bf (B)}.
\end{description}
\end{theorem} 

\end{en-text}
%%%%%%%%%%%%%%%%%%%%%%%%%%%%%%%%%%%%%%%%%%%%%%%%%%%%%%%%%%%%%%%%%%%%%

\section{Numerical studies}\label{201229-1}
In this section we run some simulation experiments to asses the quality of the estimator of the change point and of the volatilities, under two different models.
We first consider the following diffusion model without drift
\begin{equation} 
X_t=
\Bigg\{
\begin{array}{ll}
X_0+\int_0^t (1+X_s^2)^{\theta_0^*} d W_s
& \mbox{ for } t\in[0,t^*)
\y
X_{t^*}+\int_{t^*}^t (1+X_s^2)^{\theta_1^*} d W_s
& \mbox{ for } t\in[{t^*},T]. 
\end{array}
\label{eq:mod1}
\end{equation}
where $t^*$ is the true change point assumed to be  $t^* = 0.6$. The true value of the parameters are
 $\theta_0^* = 0.2$ and $\theta_1^* = \theta_0^*+n^{-\gamma}$, with $\gamma=\frac14$, $n$ is the sample size and $T=nh=1$. 
%Hence $\theta_n = |\theta_0 - \theta_1|$. 
The initial value is $X_0$ assumed to be constant, 
in particular we take $X_0 = 5$.
The sequences $a_n = b_n = \frac{1}{n \theta_n^\delta}$ with $\delta = 3$ so that they satisfy the properties required in Section \ref{081227-1}.
The first stage estimator of $\theta_0^*$ (resp. $\theta_1^*$) is obtained using the first $n a_n$ observations from the left (resp. $n a_n$ from the right). We denote the first stage estimators with $\hat \theta_i$, $i=0,1$.
Once the first stage estimators of $\theta_0^*$ and $\theta_1^*$ are available, the first stage estimator of $t^*$, i.e. $\hat t_n$ is obtained via
\beas 
\Phi_n(\hat{t}_n;\hat{\theta}_0,\hat{\theta}_1)
&=&
\min_{t\in[0,T]} 
\Phi_n(t;\hat{\theta}_0,\hat{\theta}_1). 
\eeas 
 Then, with the first stage estimator of $t^*$ in hands, we calculate the second stage estimator of $\theta_i^*$ using observations in the interval $[0, \hat t_n-b_n]$ for $\theta_0^*$ and observations in the interval $[\hat t_n + b_n, T]$ for $\theta_1^*$. 
We denote the second stage estimators of $\theta_i^*$ by $\check\theta_i$. Finally, the second stage estimator of $t^*$, i.e. $\check t_n$, is obtained as
\beas 
\Phi_n(\check{t}_n;\check{\theta}_0,\check{\theta}_1)
&=&
\min_{t\in[0,T]} 
\Phi_n(t;\check{\theta}_0,\check{\theta}_1). 
\eeas 

For comparison, we also report the value of the estimator $\tilde t_n$ obtained plugging the true parameter values in the contrast function, i.e. when the volatilities are supposed to be known
\beas 
\Phi_n(\tilde t_n;\theta_0^*,\theta_1^*)
&=&
\min_{t\in[0,T]} 
\Phi_n(t;\theta_0^*,\theta_1^*), 
\eeas 
and this can be considered as a benchmark. 
For the Monte Carlo setup, we consider different sample sizes $n=1000, 2000, 5000$ and for each sample size $n$, we run $M= 10000$ Monte Carlo replications. Under this choice of $n$
 the value of $\theta_1^* = 0.3778$, $0.3495$, and $0.3189$ respectively.
 The values of the sequences $a_n$ and $b_n$ are reported in Table \ref{tab1}.
Observations are supposed to be sampled at sample rate  $h = 1/n$.
Table \ref{tab1} also reports Monte Carlo estimates (i.e. average over the $M$ replications) of the volatility parameters $\theta_0$ and $\theta_1$ and the change point $t^*$. 
In parenthesis {\colord{are}} 
the standard deviations of the Monte Carlo estimates.
In the second experiment we consider a Cox-Ingersoll-Ross (1985) model
\begin{equation}
X_t=
\Bigg\{
\begin{array}{ll}
X_0+\int_0^t \sqrt{\theta_0^* X_s} d W_s
& \mbox{ for } t\in[0,t^*)
\y
X_{t^*}+\int_{t^*}^t \sqrt{\theta_1^* X_s} d W_s
& \mbox{ for } t\in[{t^*},T]. 
\end{array}
\label{eq:mod2}
\end{equation}
 with change point $t^*=0.7$ and all remaining experimental conditions are the same as in previous experiment. The results are reported in Table \ref{tab2}. The difference in the two experiments is only in the regularity of the diffusion coefficient term.
Comparing the two simulation results, it is possible to see that the second stage estimators in the second experiment performs slightly better in term of the standard deviation.

\begin{figure}
\centering{\includegraphics[width=0.9\textwidth,height=0.4\textheight]{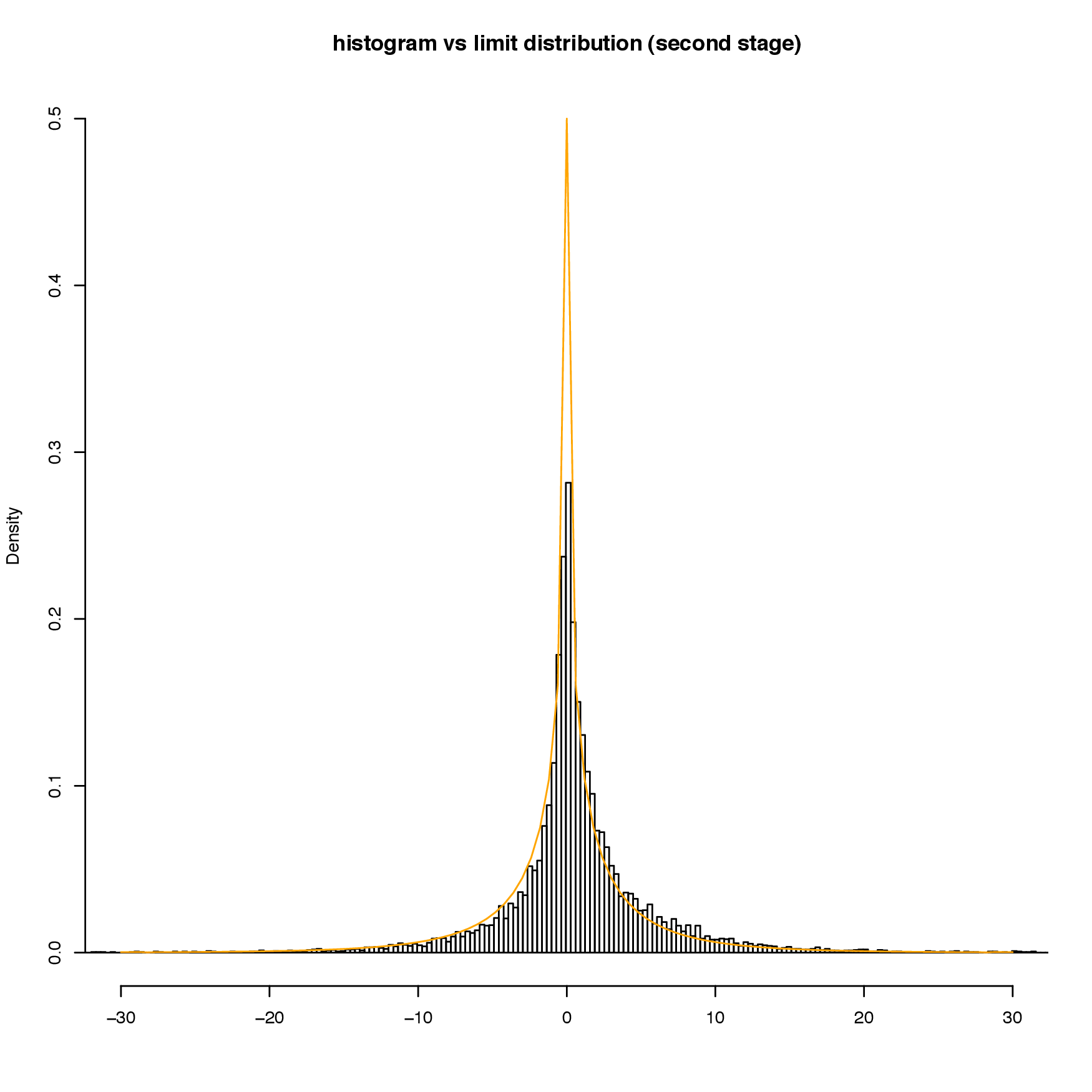}}
\centering{\includegraphics[width=0.9\textwidth,height=0.4\textheight]{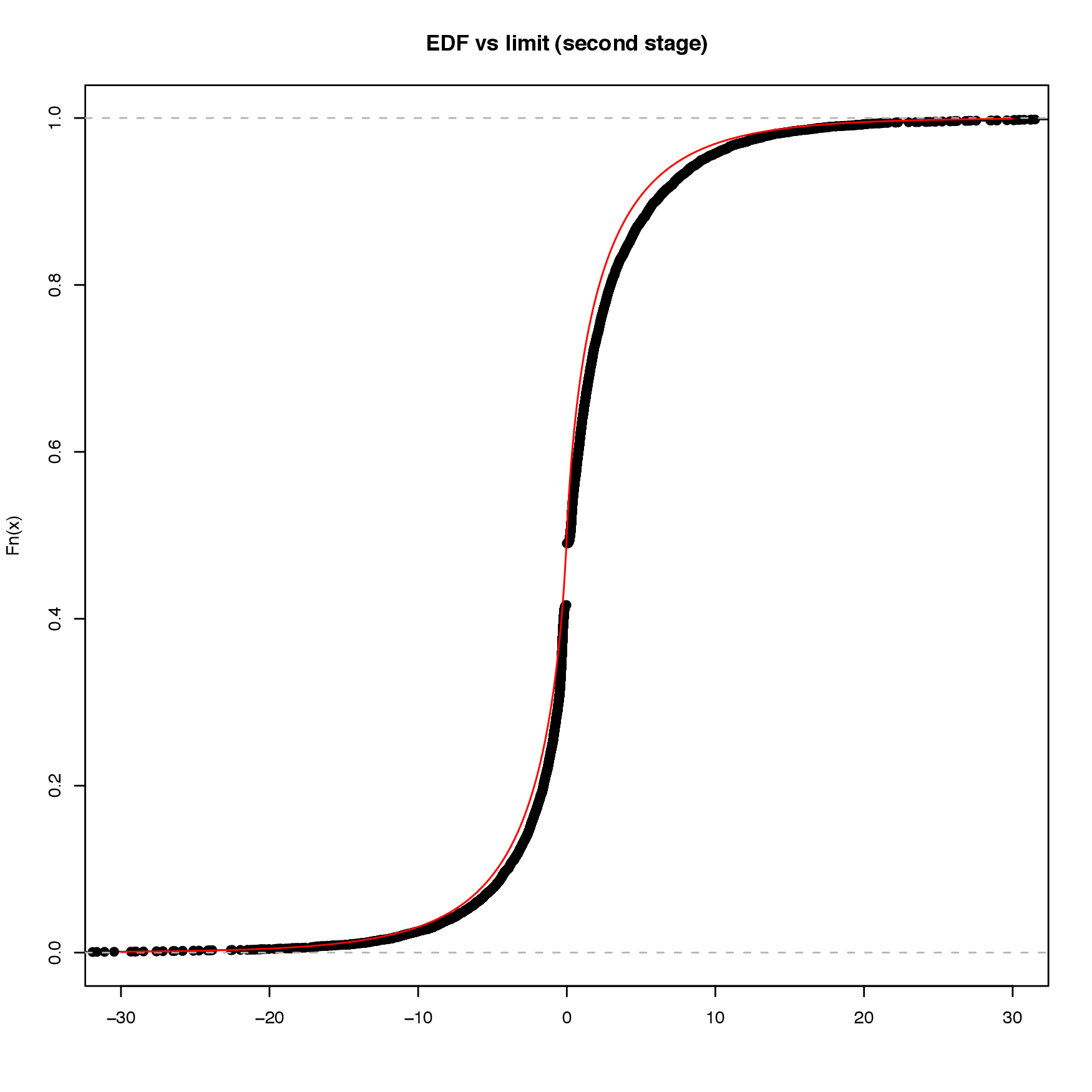}}
\caption{Histogram versus theoretical density function (up) and empirical distribution function versus theoretical distribution function (bottom) for the second stage change point estimator. Results of 10000 Monte Carlo replications and sample size $n=5000$ for the first model.}
\label{fig1}
\end{figure}

{\colorb  We also consider the behaviour of the asymptotic distribution of the change point estimator
for second stage estimator in the first model, for sample size $n=5000$.
In particular, due to mixed-normal limit, we studied the distribution of the studentized limiting distribution of  $n\theta_n^2(\check{t}_n-t^*)$ under the true model, i.e.
$$
Z = n\theta_n^2(\check{t}_n-t^*)\hat\Gamma(X_{t^*},\theta_0),
$$
with $\hat\Gamma(X_{t^*},\theta_0) = ( \log(1+X_{t^*}^2))^2$. Then $Z$ converges to
$\calw(v)-\frac12 |v|$ with density 
$$
f(x) = \frac{3}{2}e^{|x|}\left(1-\Phi\left(\frac32\sqrt{|x|}\right)\right)-\frac12\left(1-\Phi\left(\frac12\sqrt{|x|}\right)\right)
$$
and distribution function 
$$
F(x)=
\begin{cases}
g(x), &x>0 \\
1-g(-x), &x\leq 0 \\
\end{cases}
$$
with $\Phi(x)$ the distribution function of the gaussian random variable,  and
$$
g(x) = 1 + \sqrt{\frac{x}{2 \pi}} e^{-\frac{x}{8}} -\frac12 (x+5)\Phi\left(-\frac{\sqrt{x}}{2}\right)+\frac32 e^x \Phi\left(-\frac32 \sqrt{x}\right)
$$
(see e.g. Cs\"{o}rg\H{o} and Horv\'{a}th, 1997).
In Figure \ref{fig1} we report the graphical representation of the histogram and empirical distribution function of $Z$ (over 10000 Monte Carlo replications) against their theoretical counterparts which looks quite reasonable.
}
\begin{table}[H]
\begin{center}
{\small
\begin{tabular}{cc|c|ccc|ccc}
$n$ & $a_n$ & $\tilde t_n$ & $\hat\theta_0$& $\hat\theta_1$& $\hat t_n$ & $\check{\theta}_0$ & $\check{\theta}_1$ & $\check t_n$\\
\hline
5000 & 0.1189 & 0.601 & 0.200 & 0.319 & 0.601 & 0.200 & 0.319 & 0.601\\
 &&(0.005) & (0.009) & (0.014) & (0.011) & (0.005) & (0.013) & (0.012)\\
&&&&&&&\\
 2000  & 0.1495 & 0.601 & 0.200 & 0.349 & 0.601 & 0.200 & 0.349 & 0.601\\
 && (0.008) & (0.013) & (0.020) & (0.014) & (0.008) & (0.017) & (0.015)  \\
&&&&&&&\\
1000 & 0.1778 & 0.601 & 0.199 & 0.377 & 0.601 & 0.200 & 0.377 & 0.602\\
&& (0.011) & (0.017) & (0.025) & (0.019) & (0.011) & (0.026) & (0.018)\\
\end{tabular}
}
\end{center}
\caption{\small Monte Carlo estimates for model \eqref{eq:mod1} over $10000$ replications. True values: $\theta_0^* = 0.2$, $\theta_1^* = 0.378$, $0.350$, and $0.319$ for different sample sizes $n=1000$, $2000$ and $5000$. True change point $t^* = 0.6$.}
\label{tab1}
\end{table}
\begin{table}[H]
\begin{center}
{\small
\begin{tabular}{cc|c|ccc|ccc}
$n$ & $a_n$ &$\tilde t_n$ & $\hat\theta_0$& $\hat\theta_1$& $\hat t_n$ & $\check{\theta}_0$ & $\check{\theta}_1$ & $\check t_n$\\
\hline
5000 & 0.1189 &0.701& 0.200& 0.319& 0.701& 0.200& 0.319 & 0.701\\
 && (0.010) & (0.012)& (0.018) & (0.011)& (0.018)& (0.012) & (0.010)\\
 &&&&&&&\\
 2000  & 0.1495 & 0.702 & 0.200 & 0.350 & 0.701 & 0.200 & 0.350 & 0.701\\
 && (0.016) & (0.016) & (0.029) & (0.024) & (0.009) & (0.030) & (0.021)\\  
&&&&&&&\\
1000 & 0.1778 & 0.703 & 0.200 & 0.378 & 0.701 & 0.200 & 0.377 & 0.701\\
&& (0.025) & (0.021) & (0.040) & (0.038) & (0.012) & (0.056) &(0.040)\\
\end{tabular}
}
\end{center}
\caption{\small Monte Carlo estimates for model \eqref{eq:mod2} over $10000$ replications. True values: $\theta_0^* = 0.2$, $\theta_1^* = 0.378$, $0.350$, and $0.319$ for different sample sizes $n=1000$, $2000$ and $5000$. True change point $t^* = 0.7$.}
\label{tab2}
\end{table}

%\fi
%%%%%%%%%%%%%%%%%%%%%%%%%%%%%%%%%%%%%%%%%%%%%%%%%%%%%%
%%%%%%%%%%%%%%%%%%%%%%%%%%%%%%%%%%%%%%%%%%%%%%%%%%%%%%
%%%%%%%%%%%%%%%%%%%%%%%%%%%%%%%%%%%%%%%%%%%%%%%%%%%%%%
%%%%%%%%%%%%%%%%%%%%%%%%%%%%%%%%%%%%%%%%%%%%%%%%%%%%%%
%%%%%%%%%%%%%%%%%%%%%%%%%%%%%%%%%%%%%%%%%%%%%%%%%%%%%%
%%%%%%%%%%%%%%%%%%%%%%%%%%%%%%%%%%%%%%%%%%%%%%%%%%%%%%
%%%%%%%%%%%%%%%%%%%%%%%%%%%%%%%%%%%%%%%%%%%%%%%%%%%%%%
%%%%%%%%%%%%%%%%%%%%%%%%%%%%%%%%%%%%%%%%%%%%%%%%%%%%%%
%%%%%%%%%%%%%%%%%%%%%%%%%%%%%%%%%%%%%%%%%%%%%%%%%%%%%%
%%%%%%%%%%%%%%%%%%%%%%%%%%%%%%%%%%%%%%%%%%%%%%%%%%%%%%
%%%%%%%%%%%%%%%%%%%%%%%%%%%%%%%%%%%%%%%%%%%%%%%%%%%%%%

\end{document}